\documentclass[12pt]{article}
\usepackage{epsfig}
\usepackage{amsmath}
\usepackage{amssymb}
\usepackage{amscd}
\usepackage{color}
\usepackage{diagbox}

\DeclareMathOperator{\argmin}{argmin}
\newcommand{\nr}[1]{\left\Vert#1\right\Vert}

\title{Computational Analysis of a Contraction Rheometer \\
for the Grade-two Fluid Model}

\date{\today}
\author{Sara Pollock\\Department of Mathematics, University of Florida\\
\\
L. Ridgway Scott\\
University of Chicago 
}

\begin{document}
\newcommand{\Rd}{\color{red}}
\newcommand{\Bk}{\color{black}}

\newcommand\critc{\gamma}
\newcommand\notgam{{\mathcal K}}
\newcommand\Id{{\mathcal I}}
\newcommand\lemm{{\mathcal L}_m}
\newcommand\treforth{\frac{3}{4}}
\newcommand\forth{{\textstyle{\frac{1}{4}}}}
\newcommand{\intox}[1]{\int_{\Omega} #1 \, dx}
\newcommand{\cutom}{\Omega_\gamma}
\newcommand{\boint}[1]{\int_{\partial\Omega} #1 \, ds}
\newcommand{\order}[1]{\mbox{${\cal O}(#1)$}}
\newcommand{\soc}{\sigma_q}
\newcommand{\nusoc}{\sigma_{\kern-.05em q}^{\kern+.05em \prime}}
\newcommand{\nusocr}{\sigma_{\kern-.05em r}^{\kern+.05em \prime}}
\newcommand{\sigr}{\sigma_{\kern-.05em r}^{\kern+.05em \prime}}
\newcommand{\cai}{{\cal I}}
\newcommand{\caj}{{\cal J}}
\newcommand{\cak}{{\cal K}}
\newcommand{\bcf}{{\cal C}}
\newcommand{\bfC}{{\cal C}}
\newcommand{\calen}{{\cal N}}
\newcommand{\omitit}[1]{}
\newcommand{\qbox}[1]{\quad\hbox{#1}\quad}
\newcommand{\bcurl}{{\bf curl}\,}
\newcommand{\curl}{{\rm curl}\,}
\newcommand{\os}{{\hbox{\tiny\rm O}}}
\newcommand{\gs}{{\rm G}}
\newcommand{\st}{{\rm S}}
\newcommand{\ns}{{\rm N}}

\newcommand{\gss}{{\hbox{\small\rm G}}}
\newcommand{\xss}{{\hbox{\small\rm X}}}
\newcommand{\xs}{{\hbox{\tiny\rm X}}}
\newcommand{\oss}{{\hbox{\small\rm O}}}
\newcommand{\eps}{\epsilon}
\newcommand{\bxi}{{\boldsymbol \xi}}
\newcommand{\bphi}{{\boldsymbol \phi}}
\newcommand{\bpsi}{{\boldsymbol \psi}}
\newcommand{\bSigma}{{\boldsymbol \Sigma}}
\newcommand{\bchi}{{\boldsymbol \chi}}
\newcommand{\bdelta}{{\boldsymbol \delta}}
\newcommand{\bgamma}{{\boldsymbol \gamma}}
\newcommand{\bsigma}{{\boldsymbol \sigma}}
\newcommand{\bsa}{{\boldsymbol \sigma}_{\kern-.05em a}}
\newcommand{\btau}{{\boldsymbol \tau}}
\newcommand{\utVee}{{\mathbf V}_h}
\newcommand{\mybox}[1]{{\fbox{$\displaystyle #1$}}}
\newcommand{\set}[2]{\left\lbrace #1 \; \big| \; #2 \right\rbrace}
\newcommand{\beginproof}{\medskip\par\noindent{\bf Proof.\ }}
\newcommand{\beginproofof}[2]{\medskip\par\noindent{\bf Proof (of #1 \ref{#2}).\ }}
\newcommand{\proofend}{{\bf QED}\par\medskip}
\newcommand{\proofendnopar}{{\bf QED}}

\newcommand{\seedqdee}{{\cal C}} 
\newcommand{\sigmaq}{{\sigma_q}}
\newcommand{\sigmaqsq}{{\sigma_q^2}}
\newcommand{\hatseetu}{{\hat\sigma}}
\newcommand{\rhs}{{\cal T}(\nabla\uu,\lambda_2,\mu_2)}
\newcommand\rhss[1]{{\cal T}(\nabla\uu^{#1},\lambda_2,\mu_2)}
\newcommand{\ww}{{\mathbf w}}
\newcommand{\www}{{\boldsymbol \omega}}
\newcommand{\curltwo}{{\rm curl}\,}

\newcommand{\sturk}{s}
\newcommand{\sfc}{C_S} 
\newcommand{\tc}{\kern-.15em\circ\kern-.15em}
\newcommand{\upM}{{\overline M}}
\newcommand{\lowM}{{\underline M}}
\newcommand{\uptau}{\tau}
\newcommand{\lowtau}{\hat\tau}
\newcommand{\domain}{{\Omega}}
\newcommand{\bo}{{\partial\Omega}}
\newcommand\form[1]{\int_\domain #1 \,dx}
\newcommand{\curlu}{{\mathbf R}}
\newcommand{\curluo}{{\mathbf R}_O}
\newcommand{\curlug}{{\mathbf R}_G}
\newcommand{\curlv}{\widetilde{\mathbf R}}
\newcommand{\nabladot}{\nabla\kern-.15em\cdot\kern-.1em}

\newcommand{\dpdph}[1]{\frac{\partial{\phantom{#1}}}{\partial{#1}}}
\newcommand\norm[1]{\Vert\, #1 \,\Vert}
\newcommand\bignorm[1]{\big\Vert\, #1 \,\big\Vert}
\newcommand\eqnn[1]{(\ref{#1})}
\newcommand{\sdiv}{{{\nabla\cdot} \,}}
\newcommand{\nsdiv}{{{\nabla\cdot}}}
\newcommand\materderiv{\frac{{\rm D}\phantom{t}}{{\rm D}t}}
\newcommand\lowerconvderiv{\frac{{\Delta}\phantom{t}}{{\Delta}t}}
\newcommand\upperconvderiv{\frac{{\nabla}\phantom{t}}{{\nabla}t}}
\newcommand\oderiv{\frac{{d}\phantom{t}}{{d}t}}
\newcommand\ptderiv{\frac{\partial\phantom{t}}{{\partial}t}}
\newcommand\ellt{{\nabla\uu}}
\newcommand\aou{\big({\nabla\uu}+{\nabla\uu^t}\big)}
\newcommand\lplt{\left(\ellt+\ellt^t\right)}
\newcommand\half{{\textstyle{\frac{1}{2}}}}
\newcommand\threehalfs{{\textstyle{\frac{3}{2}}}}
\newcommand\third{{\textstyle{\frac{1}{3}}}}
\newcommand\ninth{{\textstyle{\frac{1}{9}}}}
\newcommand\ateninth{{\textstyle{\frac{8}{9}}}}
\newcommand\twothirds{{\textstyle{\frac{2}{3}}}}
\newcommand\fourth{{\textstyle{\frac{1}{4}}}}
\newcommand \bfz{{\bf 0}}
\newcommand \gbc{{\bf g}}
\newcommand \Rv{\hat{\bf R}}
\newcommand \Ev{{\widetilde {\bf E}}}
\newcommand{\vf}{\varphi}
\newcommand{\cn}{{\bf C}^n}
\newcommand{\intrn}{\int_{I\!\!{R^3}}}
\newcommand{\pdd}{\partial}
\newcommand {\coi}{C_0^\infty}
\newcommand {\coiof}[1]{\coi\left(#1\right)}
\newcommand {\coirn}{\coi(\rn)}
\newcommand{\ve}[3]{(#1_#2,\ldots,#1_{#3})}
\newcommand{\cbn}[2]{\left(\begin{array}{c}#1\\#2\end{array}\right)}
\newcommand \ddt{\frac{d\,}{dt}}
\newcommand{\pdx}[1]{\frac{\partial\,}{\partial #1}}
\newcommand{\ddx}[2]{\frac{\partial #1}{\partial #2}}
\newcommand{\ddsnd}[3]{\frac{\partial^2 #1}{\partial #2\partial #3}}
\newtheorem{theorem}{Theorem}[section]
\newtheorem{lemma}[theorem]{Lemma}
\newtheorem{corollary}[theorem]{Corollary}
\newenvironment{proof}{{\it Proof.}}{\hfill\qed\newline}
\newtheorem{definition}[theorem]{Definition}
\newtheorem{proposition}[theorem]{Proposition}
\newtheorem{algorithm}[theorem]{Algorithm}
\newcommand{\hsu}{{\hspace{.5cm}}}
\newcommand{\hsd}{{\hspace{1cm}}}
\newcommand{\BR}{\mathbb{R}}
\newcommand{\Rdee}{\mathbb{R}^d}
\newcommand{\Rtwo}{\mathbb{R}^2}
\newcommand{\Rtre}{\mathbb{R}^3}
\newcommand{\Rarr}{\mathbb{R}^r}
\newcommand{\Remm}{\mathbb{R}^m}
\newcommand{\ff}{{\mathbf f}}
\newcommand{\uu}{{\mathbf u}}
\newcommand{\nn}{{\mathbf n}}
\newcommand{\bt}{{\mathbf t}}
\newcommand{\ee}{{\mathbf e}}
\newcommand{\vv}{{\mathbf v}}
\newcommand{\xx}{{\mathbf x}}
\newcommand{\yy}{{\mathbf y}}
\newcommand{\zz}{{\mathbf z}}
\newcommand{\rr}{{\mathbf r}}
\newcommand{\g}{{\mathbf g}}
\newcommand{\f}{{\mathbf f}}
\newcommand{\n}{{\mathbf n}}
\newcommand{\w}{{\mathbf w}}
\newcommand{\A}{{\mathbf A}}
\newcommand{\J}{{\mathbf J}}
\newcommand{\K}{{\mathbf K}}
\newcommand{\AG}{{(\nabla\uu +\nabla\uu^t)}}
\newcommand{\B}{{\mathbf B}}
\newcommand{\G}{{\mathbf G}}
\newcommand{\Ss}{{\mathbf S}}
\newcommand{\T}{{\mathbf T}}
\newcommand{\U}{{\mathbf U}}
\newcommand{\V}{{\mathbf V}}
\newcommand{\W}{{\mathbf w}}
\newcommand{\wno}{{\tiny\rm W}}
\newcommand{\E}{{\mathbf E}}
\newcommand{\M}{{\mathbf M}}
\newcommand{\N}{{\mathbf N}}
\newcommand{\PT}{{\mathbf P}}
\newcommand{\hPT}{{\hat{\mathbf P}}}
\newcommand{\ta}{{\mathbf t}}
\newcommand{\la}{{\mathbf \lambda}}
\newcommand {\mm}{{\mathbf \mu}}
\newcommand{\cc}{{\bf curl}}
\newcommand{\dd}{\mbox{div}}
\renewcommand{\theequation}{\thesection.\arabic{equation}}
\newcommand{\strate}{{\dot{\varepsilon}}}

\newcommand{\balpha}{{\boldsymbol \alpha}}
\newcommand{\br}{{\boldsymbol r}}
\newcommand{\bdf}{{\boldsymbol f}}
\newcommand{\bU}{{\boldsymbol U}}
\newcommand{\balfo}{{\boldsymbol \alpha_1}}
\newcommand{\balft}{{\boldsymbol \alpha_2}}
\newcommand{\cI}{{\mathcal I}}

\maketitle

\begin{abstract}
We explore the possibility of simulating the
grade-two fluid model in a geometry related to a contraction rheometer,
and we provide details on several key aspects of the computation.
We show how the results can be used to determine the viscosity $\nu$
from experimental data.
We also explore the identifiability of the grade-two parameters
$\alpha_1$ and $\alpha_2$ from experimental data.
In particular, as the flow rate varies, force data appears to be nearly the same for 
certain distinct pairs of values $\alpha_1$ and $\alpha_2$; however we determine a 
regime for $\alpha_1$ and $\alpha_2$ for which the parameters may be identifiable
with a contraction rheometer.
\end{abstract}

\section{Introduction}\label{sec:intro}

A rheometer is a device that is used to determine physical properties of fluids
such as viscosity, but also other properties for non-Newtonian fluids,
via controlled experiments.
To use a rheometer to characterize parameters used in computational
models of fluids, it is necessary to have a way to convert data from
the form produced by the rheometer into an estimate of model parameters.
In ideal cases of hypothetical rheometers \cite{lrsBIBix},
this can be done with analytical solutions of the model equations.
But in more realistic cases, it is necessary to solve the model equations
numerically in the rheometer geometry, creating a mapping from model parameters to
approximations of the rheometer data.
Then one can attempt to invert this mapping to generate model parameters
from experimental data.
Here we explore the first step (the forward problem) for a particular
model and a particular rheometer geometry.
In a subsequent paper \cite{lrsBIBkq}, we will examine 
the corresponding inverse problem.

The grade-two fluid model is the lowest-order member of a family of models 
proposed by Rivlin and Ericksen \cite{ref:RivlinEricksen,lrsBIBej}.
In these models,
the stress-strain relationship involves derivatives of the fluid velocity.
The grade-two model involves two parameters in addition to the fluid viscosity
and has been widely studied \cite{VGrheoBook}.
On the other hand, computational models have been limited so far by restrictions
on the two parameters, the dimension, or boundary conditions
\cite{ref:gradetwoalphatwo,bernard2018fully,lrsBIBem,lrsBIBej}.

Recently, an algorithm was proposed \cite{lrsBIBjd} for both two and three dimensions
that supports the use of general parameters and inflow boundary conditions.
We use this algorithm to compute a force integral over the contraction region
for the flow in a contracting duct.  
We compute the force integral over a range of parameters and carefully analyze the
structure of the resulting data to determine a regime in which we may be able to 
identify parameters in the grade-two model.

We study the proposed computational method in detail.
We point out another algorithm designed for more restrictive settings.
The latter algorithm is able to find solutions for much larger parameter values,
for constrained parameter values ($\alpha_1+\alpha_2=0$).
To improve our solves and to extend the parameter regime where we are able to solve
with the more general method, we consider an implementation of the extrapolation technique
known as Anderson Acceleration (AA)
\cite{anderson65}, similar to that given in \cite{PR21,PR23}.

The remainder of this paper is organized as follows.
Sections \ref{subsec:g2model}-\ref{subsec:theord} summarize the main results from 
\cite{lrsBIBjd}, namely the new algorithm for the grade-two model and the theory
describing its convergence;
{additionally, here we provide details on an Anderson accelerated version
of the nonlinear iteration.
Section \ref{sec:cduct-flow} describes the contracting duct 
domain over which we will perform the force integral computations, and describes the main
characteristics of the flow in this duct.  
Section \ref{sec:rheometer} contains the main results of the paper.
In this section we define the force integral and provide a study
of the structure of the computed data with respect to each of its
parameters in order to determine a regime in which this rheometer
may be used to identify parameters in the grade-two model.
Section \ref{sec:compudet} contains the computational details of the calculations 
performed in the preceeding section focusing on the computational mesh and its local 
refinement. The included appendix contains further technical details on determining
appropriate inflow boundary conditions for grade-two flow in both a channel and a pipe.

\subsection{Grade-two fluid model}\label{subsec:g2model}

The lowest-order, grade-two model of Rivlin and
Ericksen \cite{ref:RivlinEricksen,lrsBIBej} can be expressed as 
\begin{equation}\label{eqn:fullgeetoo}
\begin{split}
-\nu\,\Delta\,\uu +\uu\cdot\nabla\uu +\nabla p &=\sdiv\widehat\btau,\\
\sdiv\uu=0  \quad\hbox{in}\;\Omega,  \quad \uu&=\gbc\quad\hbox{on}\;\partial\Omega ,
\end{split}
\end{equation}
where 
\begin{equation}\label{eqn:sigmahatdo}
\begin{split}
\widehat\btau
&=\alpha_1\big(\uu\cdot\nabla\A-\A\tc(\nabla\uu)^t-(\nabla\uu)\tc\A\big) 
    + (2\alpha_1+\alpha_2) \A\tc\A,
\end{split}
\end{equation}
and $\A=\nabla\uu+(\nabla\uu)^t$.
We assume that the boundary data $\gbc$ is defined on all $\Omega$, is divergence
free, and sufficiently smooth, to be specified subsequently.

The equations (\ref{eqn:fullgeetoo}--\ref{eqn:sigmahatdo}) can be viewed as a
perturbation of the Navier--Stokes system
\begin{equation}\label{eqn:navstonlyo}
\begin{split}
-\Delta\,\uu +R\uu\cdot\nabla\uu +\nabla p &=\bfz,\\
\sdiv\uu=0  \quad\hbox{in}\;\Omega,  \quad \uu&=\gbc\quad\hbox{on}\;\partial\Omega ,
\end{split}
\end{equation}
where $R=UL/\nu$ is the Reynolds number ($U$ is a velocity  scale and $L$ is a
length scale used to nondimensionalize the equations).
When $R=0$, the system \eqref{eqn:navstonlyo} is called the Stokes equations.
The pressure $p$ has been rescaled as well.

\subsection{Special case in dimension two}\label{subsec:special}

When $\alpha_1+\alpha_2=0$, there is a simplification \cite{lrsBIBem,lrsBIBej}
that can be made in two dimensions that reduces the system 
(\ref{eqn:fullgeetoo}--\ref{eqn:sigmahatdo}) to
\begin{equation}\label{eqn:simpgtoodeeo}
\begin{split}
-\nu\Delta\uu+z(u_2,-u_1)^t+\nabla q = \bfz, \quad \sdiv\uu=0,\quad 
\nu z+\alpha_1\uu\cdot\nabla z=\nu\,\curl\uu,\quad \hbox{in}\;\Omega ,
\end{split}
\end{equation}
where $\curl\uu=(u_{2,1}-u_{1,2})$ and $p=q+\half|\uu|^2$.
Unfortunately, this simplification does not generalize to $\alpha_2\neq - \alpha_1$
or to three dimensions.
Thus we consider the following.

\subsection{Alternate formulation of the grade-two model equations}

Define the tensor $\uu\otimes\uu$ by $(\uu\otimes\uu)_{ij}=u_i\,u_j$.
Then
\begin{equation*}
\sdiv(\uu\otimes\uu)=\uu\cdot\nabla\uu.
\end{equation*}

Let $\pi$ be related to $p$ by
\begin{equation}\label{eqn:peetopigoo}
\nu\pi+\alpha_1 \uu\cdot\nabla\pi=p.
\end{equation}
Define
\begin{equation}\label{eqn:sigmadefoo}
\btau = \alpha_1(\nabla\uu)^t\tc\A +(\alpha_1+ \alpha_2)\A\tc\A-\uu\otimes\uu,
\end{equation}
and
\begin{equation*}
N(\uu,\pi)= -\alpha_1\pi\nabla\uu^t + \btau .
\end{equation*}
Note that $N$ is not a symmetric tensor due to the term $\pi\nabla\uu^t$.
The incompressibility condition $\sdiv\uu=0$ implies that
\begin{equation}\label{eqn:divennetoo}
\sdiv(\pi\nabla\uu^t) = \nabla\uu^t\nabla\pi,\qquad
\sdiv N(\uu,\pi)= -\alpha_1\nabla\uu^t\nabla\pi +\sdiv\btau.
\end{equation}
Therefore
\begin{equation}\label{eqn:aldiventoo}
\sdiv N(\uu,\pi)= -\alpha_1\nabla\uu^t\nabla\pi +\sdiv
\big(\alpha_1 \nabla\uu^t\tc\A +(\alpha_1+ \alpha_2) \A\tc\A   -\uu\otimes\uu\big). \\
\end{equation}

Now consider the problem proposed in \cite{lrsBIBjd}:
\begin{equation}\label{eqn:diffrfsgeetoo}
\begin{split}
-\Delta\uu &+\nabla\pi =\ww\quad\hbox{in}\;\Omega, \quad
\sdiv\uu=0  \quad\hbox{in}\;\Omega,  \quad \uu=\gbc\quad\hbox{on}\;\partial\Omega, \\
(\nu I&+\alpha_1 \uu\cdot\nabla)\ww =\sdiv N(\uu,\pi) \quad\hbox{in}\;\Omega,\quad
\ww=\ww_b\quad\hbox{on}\;\Gamma_-,
\end{split}
\end{equation}
where 
\begin{equation}\label{eqn:defgaminsq}
\Gamma_{-}=\set{\xx\in\partial\Omega}{\alpha_1\,\gbc(\xx)\cdot\nn<0}.
\end{equation}
The following is proved in \cite{lrsBIBjd}.

\begin{theorem}\label{thm:diffrfsgeetoo}
Suppose that $(\uu,\pi)$ solves \eqref{eqn:diffrfsgeetoo} 
and $p$ is given by \eqref{eqn:peetopigoo}.
Then $(\uu,p)$ satisfies
\eqref{eqn:fullgeetoo} with $\widehat\btau$ defined by \eqref{eqn:sigmahatdo}.
The vector function $\ww$ satisfies
$$
\ww=\frac{1}{\nu}\big(\sdiv\widehat\btau - \uu\cdot\nabla\uu - \nabla p \big) +\nabla\pi.
$$
\end{theorem}

One modeling challenge arises because \eqref{eqn:fullgeetoo} is a third-order PDE due
to the presence of the term $\uu\cdot\nabla(\Delta\uu)$.
Therefore it is necessary to specify an additonal boundary condition beyond what
would be done for the Navier--Stokes equations to have a unique solution.
The quantity $\ww$ on which we pose a boundary condition is the divergence of the stress.

\subsection{An algorithm for the transformed equations}\label{subsec:disc-alg}

The system \eqref{eqn:diffrfsgeetoo} is analogous to the reduced system
in \cite{lrsBIBej}, and that paper suggested the algorithm used
in \cite{lrsBIBjd} for solving \eqref{eqn:diffrfsgeetoo}: 
start with some $\ww^0$, then solve for $n\geq 1$
\begin{equation}\label{eqn:algfrfsgeetoo}
\begin{split}
-\Delta\uu^n &+\nabla\pi^n =\ww^{n-1}\quad\hbox{in}\;\Omega,\quad
\sdiv\uu^n=0 \quad\hbox{in}\;\Omega,  \quad \uu^n=\gbc\quad\hbox{on}\;\partial\Omega,\\
(\nu I&+\alpha_1 \uu^n\cdot\nabla)\ww^n =\sdiv N(\uu^n,\pi^n) \quad\hbox{in}\;\Omega,
\quad \ww^n=\ww_b\quad\hbox{on}\;\Gamma_-.
\end{split}
\end{equation}
For definiteness, we will take $\ww^0=\ww_b$.
In \cite{lrsBIBej}, convergence of this iteration is proved for small data $\gbc$ and $\ww_b$.

The following discrete variational form of \eqref{eqn:algfrfsgeetoo} 
which is suitable for finite element approximation, and accompanying numerical
algorithm, is given in \cite{lrsBIBjd}.
Let $W_h$ be the space of continuous, vector-valued, piecewise polynomials
of degree $k$, let $V_h=\set{\vv\in W_h}{\vv=\bfz\;\hbox{on}\;\partial\Omega}$, and let
$\Pi_h$ be continuous, scalar-valued, piecewise polynomials of degree $k-1$.
For the computations shown throughout this paper, we use $k=4$.

First, using the iterated penalty method:
find $\uu^{n,\ell}\in V_h+\gbc$ such that
\begin{equation}\label{eqn:varuipm}
\begin{split}
\intox{\nabla\uu^{n,\ell}:\nabla\vv} + \rho \intox{\sdiv\uu^{n,\ell}\,\sdiv\vv}&=
\intox{\ww^{n-1}\cdot\vv}
-\intox{\sdiv\zz^\ell\,\sdiv\vv} \quad\forall \vv\in V_h, \\
\zz^{\ell+1}&=\zz^\ell+\rho\,\uu^{n,\ell} . \\
\end{split}
\end{equation}
Once this is converged, we set $\uu^{n}=\uu^{n,\ell}$ and define the pressure via
\cite{lrsBIBia}
\begin{equation}\label{eqn:varpress}
\intox{\pi^n \, q}=\intox{-\sdiv\zz^{\ell+1} \, q} \quad\forall q\in \Pi_h.
\end{equation}

We can pose the transport equation \eqref{eqn:algfrfsgeetoo} via:
find $\ww^{n}\in \widetilde{V}_h+\ww_b$ such that
\begin{equation}\label{eqn:varw}
\begin{split}
\nu \intox{\ww^n\cdot\vv}+\alpha_1\intox{\big(\uu^n\cdot\nabla\ww^n\big)\cdot\vv}
-\intox{\big(\sdiv N(\uu^n,\pi^n)\big) \cdot\vv}=0\quad\forall \vv\in \widetilde{V}_h,
\end{split}
\end{equation}
where $\ww_b$ is posed only on the inflow boundary, that is,
$$
\widetilde{V}_h=\set{\vv\in W_h}{\vv=\bfz\;\hbox{on}\; \Gamma_-},
\qquad \Gamma_- =\set{\xx\in \partial\Omega}{\nn\cdot\gbc <0}.
$$
\subsection{Anderson accelerating the solution sequence}

We consider augmenting the solver for nonlinear system \eqref{eqn:algfrfsgeetoo},
which is implemented via \eqref{eqn:varuipm}-\eqref{eqn:varw},
by applying a filtered version of AA as in \cite{PR21,PR23}
to the approximation sequences $\{\uu^n\}$ and $\{ \zz^n \}$.

The first equation in \eqref{eqn:algfrfsgeetoo} is solved by the iterated penalty
method \eqref{eqn:varuipm}, which also generates $\zz^n$. The $L_2$ projection of 
the divergence of $\zz^n$ is computed via \eqref{eqn:varpress} and used to solve
for the auxiliary variable $\ww^n$ in \eqref{eqn:varw}. 
For consistency, 
it makes sense to perform the extrapolation on $\zz^n$ along with $\uu^n$, with 
the extrapolation parameter determined entirely by $\uu^n$.  

Specifically, we consider the following modification to 
\eqref{eqn:varuipm}-\eqref{eqn:varw}.
Starting with an initial iterate $\uu^0$, 
let $\widehat \uu^n = \uu^{n,l}$ upon convergence of \eqref{eqn:varuipm}.  
The algorithm recombines up to $m_{max}$ previous iterates
$\widehat \uu^j$ and 
update steps $\delta \uu^j = \widehat \uu^j - \uu^{j-1}$, for $n-m \le j \le n$,
to form the next iterate $\uu^n$.
It also recombines the corresponding iterates and updates steps
$\zz^j$ and $\delta \zz^j = \widehat \zz^j - \zz^{j-1}$ to form the next iterate $\zz^n$.
The algorithm with depth $m=0$ reduces to the original fixed-point iteration without
acceleration.

The filtered version of AA described below was introduced in \cite{PR21} and built 
upon in \cite{PR23} to better 
control the accumulation of higher-order terms in the residual expansion by enforcing
a sufficent linear indepenedence condition (or, if the parameter $\sigma$ is chosen
close enough to 1, a near-orthogonality condition), between the columns of the 
coefficient matrix of the underlying least-squares problem.
The acceleration both controls the growth of the iteration
count for smaller parameter pairs $(\alpha_1, \alpha_2)$ as the mesh is 
refined, and it
becomes an enabling technology allowing the solution for larger parameter pairs
$(\alpha_1,\alpha_2)$ than can be solved without the acceleration.  
The filtering technique is seen both
to decrease sensitivity to choice of the 
extrapolation depth $m$ and to reduce the number of iterations in the solve on 
finer meshes.

Filtered AA is implemented and described naturally 
as a linear algebra routine, operating on the
coefficients $U^n$ of the basis expansion 
$\uu^n = \sum U^n_i \varphi_i$,
where the $\{\varphi_i\}$ span
the discrete space $V_h$. 
In agreement with standard practice, the inner optimization for this problem is 
performed with respect to the $l_2$ norm. 

\begin{algorithm}\label{alg:FAA}
{\bf(Filtered AA.)} Set depth $m_{max}$.
Compute $\widehat U^1$ and $\delta U^n = \widehat U^1  - U^0$. 
\\
Set $m_0 = 0$, 
$F_0  = \begin{pmatrix}(\delta U^{n+1}-\delta U^n) \end{pmatrix}$ and 
$E_0 = \begin{pmatrix}(U^n - U^{n-1})\end{pmatrix}$.
\\ \noindent
For $n = 1, 2, \ldots$, set $m_n = \min\{m_{n-1}+1,m_{max}\}$
\\ \indent
Compute $\widehat X^{n+1}$ and $\delta X^{n+1} = \widehat X^{n+1}  - X^{n}$, for 
$X = \{U,Z\}$
\\ \indent
Set 
$FX_n  = \begin{pmatrix}(\delta X^{n+1}-\delta X^n)&  FX_{n-1} \end{pmatrix}$ and 
$EX_n = \begin{pmatrix}(X^n - X^{n-1}) & EX_{n-1}\end{pmatrix}$, for $X = \{U,Z\}$.
\\ \indent
Set $(EU_n,FU_n, EZ_n,FZ_n,m_n,\gamma_n)$ = 
Filter $(EU_n,FU_n,\delta U^{n+1},EZ_n, FZ_n,m_n)$.
\\ \indent
Set damping factor $0 < \beta_n \le 1$.
\\ \indent
Set $X^{n+1} = X^n + \beta_n \, \delta X^{n+1} - \left(EX_n + \beta_n \, FX_n \right)
\gamma_{n}$, for $X = \{U,Z\}$.
\end{algorithm}

The filtering algorithm computes the solution to a least-squares problem of the form
$\gamma_{n} = \argmin_{\gamma \in \BR^{m_n}} \nr{ \delta U_{n+1} - F_n \gamma}_{l_2}$,
such that columns of $F_n$ are filtered out if the direction sine between any column 
of $F_n$ and the subspace spanned by the columns to its left are
less in magnitude than parameter $\sigma$.
In \cite{PR23} this is referred to as {\em angle filtering}.
Setting $\sigma = 0$ means no filtering is performed, and setting $\sigma =1$ 
filters out any column of $F_n$ that is not orthogonal to the columns to its left.
Here we use a dynamic filtering strategy as was shown effective in 
\cite{PR21,PR23}. This method starts with a higher filtering tolerance $\sigma_{max}$ 
which filters out more columns in the preasymptotic regime and relaxes to a lower 
tolerance $\sigma_{min}$ which uses more columns for a better optimization in the 
asymptotic regime.  Here we use $\sigma_{min} = 0.1$, $\sigma_{max} = 2^{-1/2}$, and
$\sigma = \max\{ \min \{ \sigma_{max}, \|\delta U^{n+1} \|_{l_2}^{1/2} \}, \sigma_{min}\}$.

\begin{algorithm}{\bf($(E,F,EZ,FZ,m,\gamma)$ = Filter($E,F,\delta U,EZ, FZ,m$).)} 
Given  minimum and maximum filtering thresholds 
$0 \le \sigma_{min} < \sigma_{max} < 1$
\\
Compute $F = QR$, the thin QR decomposition of $F$
\\ 
Set
$\sigma = \max\{ \min \{ \sigma_{max}, \|\delta U^{n+1} \|_{l_2}^{1/2} \}, \sigma_{min}\}$
\\ 
For $i = 2, \ldots m_n$
\\ \indent
Compute $\sigma_i = |r_{ii}|/\|f_i\|_{l_2}$, 
where $r_{ii}$ is the
diagonal entry of $R$, and $f_i$ is column $i$ of $F$
\\ \indent 
If $\sigma_i < \sigma$, remove column $i$ from $F$, $E$, $FZ$ and $EZ$, 
and set $m = m-1$
\\
If any columns were removed, recompute $F = QR$
\\
Solve $R \gamma = Q^T \delta U$ for $\gamma$
\end{algorithm}

In the examples that follow we perform the iterations without damping ($\beta_n = 1$ for
all $n$).  In practice, $0 < \beta_n < 1$ can often be used to solve problems for 
a wider range of parameters.

{
\subsection{Required inflow boundary conditions}

One feature of the proposed method \eqref{eqn:algfrfsgeetoo} is that it clarifies
the required additional boundary condition, namely for $\ww=-\Delta\uu+\nabla\pi$.
Although we cannot say how to pick this in general, Appendix \ref{sec:detinflo}
computes $\ww$ for typical flow geometries.
We can extend this using Amick's theorem \cite{ref:amicksteadypipes} as described
in \cite{lrsBIBjk}.

\subsection{Theoretical details}
\label{subsec:theord}

We collect in Appendix \ref{sec:spaces}
details on the Lebesgue and Sobolev spaces and norms used.
Let $d$ be the dimension of $\Omega$.
Assume that the domain regularity asssumtion \cite[(4.2)]{lrsBIBjd} holds
for $Q_0>d$ and $Q_1>d/2$, as follows.
Suppose that the solution of
\begin{equation}\label{eqn:emptyeetoo}
-\Delta\,\uu +\nabla p =\ww \quad\hbox{and}\quad
\sdiv\uu=0 \quad\hbox{in}\;\Omega,\quad \uu=\gbc\quad\hbox{on}\;\partial\Omega,
\end{equation}
satisfies, for $1\leq q\leq Q_s$, $s=0,1$, and any $\gbc\in W^{s+2}_q(\Omega)$
and $\ww\in W^{s}_q(\Omega)$, the following estimate:
\begin{equation}\label{eqn:regasuetoo}
\norm{\uu}_{W^{s+2}_q(\Omega)}+
\norm{\pi}_{W^{s+1}_q(\Omega)}
\leq c_{q,s}\big(\norm{\ww}_{W^{s}_q(\Omega)} +
\norm{\gbc}_{W^{s+2}_q(\Omega)} \big),
\end{equation}
for a constant $c_{q,s}$ that depends only on $q$ and $s$.
This assumption holds if we round off the corners of the contracting duct.
The following is proved in \cite[Theorem 4.1]{lrsBIBjd}.

\begin{theorem}\label{thm:bigtheo}
Suppose that $d<q<Q_0$.
If the boundary data and initial iterates are sufficiently small,
the iterates \eqref{eqn:algfrfsgeetoo} are bounded for all $n>0$:
\begin{equation}\label{eqn:bigtheo}
\norm{\ww^{n}}_{L^q(\Omega)}\leq \notgam,\qquad
\norm{\uu^{n}}_{W^{2}_q(\Omega)}+
\norm{\pi^{n}}_{W^{1}_q(\Omega)}\leq c_q\big(\norm{\gbc}_{W^{2}_q(\Omega)} +\notgam\big),
\end{equation}
where $\notgam$ is a finite positive constant.
Suppose further that $r\leq Q_1$ satisfies
\begin{equation}\label{eqn:qrfacto}
\frac2d > \frac1r > \frac1q+\frac12.
\end{equation}
Then
\begin{equation}\label{eqn:medtheo}
\norm{\ww^{n}}_{W^1_r(\Omega)}\leq \notgam,\qquad
\norm{\uu^{n}}_{W^{3}_r(\Omega)}+
\norm{\pi^{n}}_{W^{2}_r(\Omega)}\leq c_q\big(\norm{\gbc}_{W^3_r(\Omega)} +\notgam\big).
\end{equation}
Moreover, $(\uu^{n},\pi^{n},\ww^{n})$ converge geometrically in 
$W^2_r(\Omega)^d\times W^1_r(\Omega)\times L^r(\Omega)^d$
to the solution $(\uu,\pi,\ww)$ of \eqref{eqn:diffrfsgeetoo}.
In view of Theorem \ref{thm:diffrfsgeetoo},
$(\uu,p)$ is the solution of the grade-two model \eqref{eqn:fullgeetoo},
where $p$ is related to $\pi$ by \eqref{eqn:peetopigoo}.
\end{theorem}

Note that there is a typo in \cite[Theorem 4.1]{lrsBIBjd},
where the estimates for $s=1$ should have $q$ replaced by $r$.

The constraint \eqref{eqn:qrfacto} implies $q>2$ for $d=2$ and $q>6$ for $d=3$,
and thus the constraint $q>d$ is satisfied implicitly.

\begin{figure}
\centerline{(a)\includegraphics[width=2.7in]{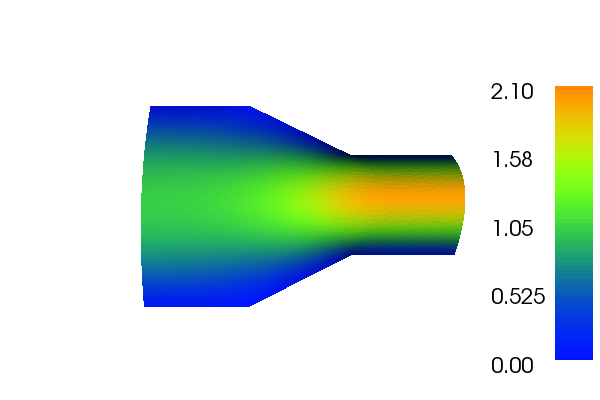}
\qquad\qquad(b)\includegraphics[width=2.7in]{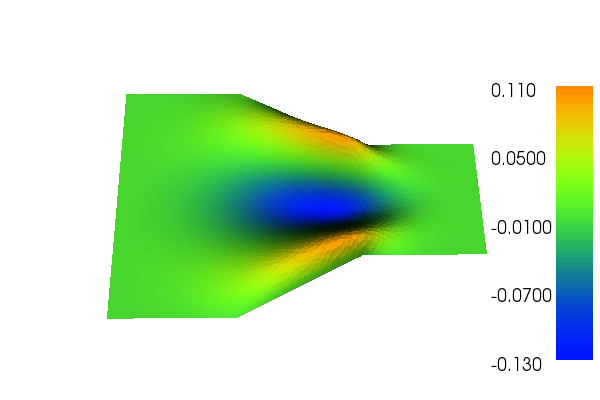}}
\caption{
Flow in a contracting duct 
{ for Stokes flow $\uu_\st$ and Navier-Stokes flow $\uu_\ns$.}
(a) horizontal flow component of $\uu_\st$, and
(b) horizontal flow component of $\uu_\ns-\uu_\st$, for
$R=10$,  both with mesh parameter 64.
The computational domain $\domain$ is as specified in \eqref{eqn:expandom}, 
with $b_i=1$, $b_o=1$, $L=1$, $H=0.5$.
Computed using \eqref{eqn:navstonlyo}.}
\label{fig:stonavstok}
\end{figure}

\section{Grade-two flow in a contracting duct}\label{sec:cduct-flow}

We begin by describing a typical flow problem involving a contracting duct.
We pose a Poiseuille flow profile at the inlet and exit of the channel.
We allow ``buffers'' at each end of the contraction for the flow to regain
the Poiseuille flow profile.
Thus the domain consists of three parts, first the inlet buffer
$$
\Omega_i=\set{(x,y)}{-b_i\leq x\leq 0,\; |y|\leq 1}.
$$
The contraction zone has length $L$ and height $H$ and is defined by
$$
\Omega_e=\set{(x,y)}{0\leq x\leq L,\; |y|\leq 1+((H-1)/L)x}.
$$
Finally, the outlet buffer is
$$
\Omega_o=\set{(x,y)}{L \leq x\leq L+b_o,\; |y|\leq H}.
$$
Then the computational domain is
\begin{equation}\label{eqn:expandom}
\domain=\Omega_i \cup \Omega_e \cup \Omega_o .
\end{equation}
We will see that the lengths of these buffer zones influence the results
substantially in some cases.

The Poiseuille-like boundary conditions we choose are as follows.
At the inlet, we choose
$$
\uu(-b_i,y)=(1-y^2,0)^t,\quad y\in[-1,1].
$$
At the outlet, we choose
$$
\uu(L+b_o,y)=\big(H^{-1}(1-(y/H)^2),0\big)^t,\quad y\in[-H,H].
$$
The corresponding inflow boundary conditions for $\ww$ \cite{lrsBIBjd} are
\begin{equation}\label{eqn:twobesee}
\begin{split}
\ww(-b_i,y) = \big(0, \frac{4U^2}{\nu}y(3 \alpha_1 + 2\alpha_2) \big)^t, \quad y \in [-1,1]. 
\end{split}
\end{equation}

\subsection{Stokes versus Navier--Stokes}

In figure \ref{fig:stonavstok}, we see the horizontal flow component 
in the domain $\domain$ for Stokes flow $\uu_\st$, shown in panel (a), and the 
horizontal flow component of the difference $\uu_\ns-\uu_\st$, shown in panel (b),
between the horizontal flow component of the Navier-Stokes solution 
($R=10$) and the Stokes solution.
Several features are of interest.
First of all, the Navier-Stokes solution returns to the parabolic
profile quickly both before and after the contraction.
Secondly, there is a significant boundary layer for the Navier-Stokes solution
in the contraction zone, and the flow there is more plug-like, with the
Stokes solution being larger in the middle of the contraction zone.
This may be counter-intuitive, in that the Stokes flow is faster in
the middle of the contraction zone, but this is consistent with what
is known for Jeffrey-Hamel flow \cite{LandauLifshitz}.

\begin{figure}
\centerline{(a)\includegraphics[width=2.7in]{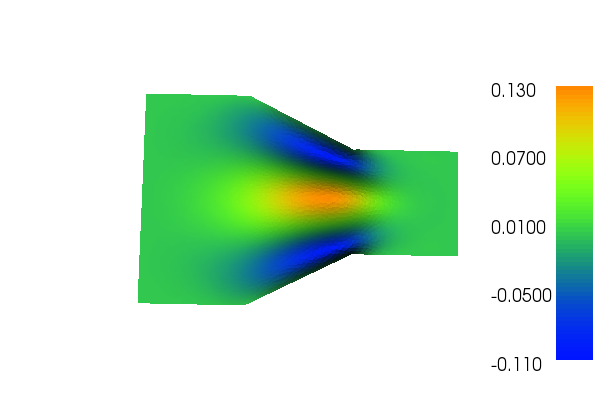}
\qquad     (b)\includegraphics[width=2.7in]{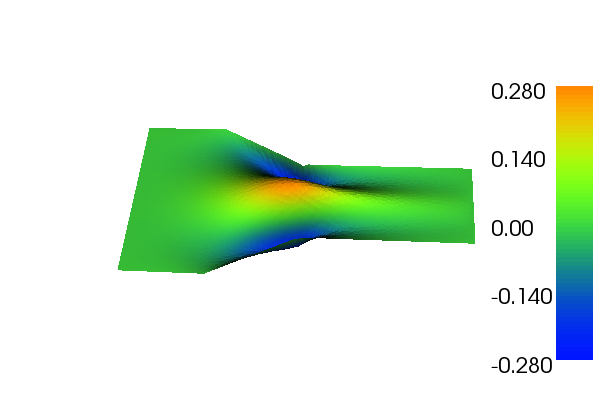}}
\caption{
Horizontal flow  of the difference $\uu_\gs-\uu_\ns$ for $\uu_\gs$
being the solution of the grade-two model \eqref{eqn:simpgtoodeeo} with
(a) $R=10$,  $\alpha_1=10$, $\alpha_2=-10$, and
(b) $R=40$, $\alpha_1=1$,  $\alpha_2=-1$, 
both with mesh parameter 64.
The computational domain $\domain$ is as specified in \eqref{eqn:expandom}, 
with $b_i=1$, $L=1$, $H=0.5$, and with (a) $b_o=1$,  (b) $b_o=2$.
Computed using \eqref{eqn:simpgtoodeeo}.}
\label{fig:navstok}
\end{figure}

\subsection{Grade-two with $\alpha_1+\alpha_2=0$}

In figure \ref{fig:navstok}(a), we depict the grade-two flow with $\alpha_2=-\alpha_1$ 
computed via \eqref{eqn:simpgtoodeeo}
in the domain $\domain$ defined in \eqref{eqn:expandom}, showing the horizontal 
flow component for the difference $\uu_\gs-\uu_\ns$ for $R=10$ and $\alpha_1=10$.
For that domain, we have $b_i=1$, $b_o=1$, $L=1$, $H=0.5$.
Here we use the shorthand Re for the Reynolds number, and for the numerical
value we write $R$.
In figure \ref{fig:navstok}(b), we depict the grade-two flow in the 
domain $\domain$ defined in \eqref{eqn:expandom}, showing the horizontal 
flow component for the difference $\uu_\gs-\uu_\ns$ for $R=40$ and $\alpha_1=1$.
For that domain, we have changed to $b_o=2$ to allow the flow to return
to a parabolic form in the outflow buffer.

\begin{table}
\begin{center} 
\begin{tabular}{|r|r|c|c|c|c|c|c|}
\hline
Re\; &$\alpha_1$\;\, &$\uu_\ns-\uu_\gs$ &$\uu_\st-\uu_\gs$ &$\uu_\ns-\uu_\st$  
&$\norm{z}_{L^2}$ & $\Delta p_\ns$ &  $\Delta p_\gs$\\
\hline
  0.1 & 0.01 & 2.26e-06 & 8.17e-04 & 8.18e-04 & 8.7577 & 575.2 & 575.3 \\
  1.0 & 0.01 & 2.25e-04 & 7.98e-03 & 8.09e-03 & 8.7497 & 58.45 & 58.55 \\
  1.0 &  0.1 & 1.87e-03 & 7.00e-03 & 8.09e-03 & 8.6597 & 58.45 & 59.45 \\
  1.0 &  1.0 & 6.07e-03 & 3.34e-03 & 8.09e-03 & 7.6517 & 58.45 & 68.49 \\
 10.0 &  1.0 & 6.64e-02 & 6.61e-03 & 6.97e-02 & 4.7639 & 6.801 & 16.83 \\
 50.0 &  1.0 & 2.09e-01 & 7.31e-03 & 2.14e-01 & 3.5618 & 2.232 & 12.24 \\
 10.0 & 10.0 & 6.93e-02 & 7.48e-04 & 6.97e-02 & 3.2948 & 6.801 & 107.2 \\
 50.0 & 10.0 & 2.14e-01 & 8.12e-04 & 2.14e-01 & 3.0423 & 2.232 & 102.6 \\
\hline
\end{tabular}
\end{center}
\caption{Relative differences 
$\norm{\uu_a-\uu_b}_{H^1(\domain)}/\norm{\uu_\st}_{H^1(\domain)}$,
indicated in columns 3--5 by $\uu_a-\uu_b$, between 
solutions to Grade-Two $\uu_\gs$, Navier-Stokes $\uu_\ns$, and Stokes $\uu_\st$.
$\norm{\uu_\st}_{H^1(\domain)}=9.2616 $ in all cases.
BC's indicates the boundary conditions on the $z$ equation.
The parameters defining the computational domain \eqref{eqn:expandom} are,
in all cases, $H=0.5$, $L=1$, $b_o=3$, and $b_i=1$, and the meshsize is 64.
In all cases, $\alpha_2=-\alpha_1$. Computed using \eqref{eqn:simpgtoodeeo}.
}
\label{tabl:plusminust}
\end{table}

Table \ref{tabl:plusminust} gives data for other values of $\alpha_1$,
but still with $\alpha_2=-\alpha_1$.
Several things emerge from this table.
First of all, it is evident that the discrepancy between the Stokes
and Navier-Stokes equations is close to linear for small Re, with
the coefficient in this case being on the order of $0.008$.
But for larger Re, the relationship is sublinear.
Similarly, the discrepancy between Navier-Stokes and grade-two is
close to linear for small $\alpha_1$, with
the coefficient in this case being on the order of $0.02$ for $R=1$.
However, for $R=10$, the difference between $\alpha_1=1$ and $\alpha_1=10$
is minimal.
The same thing is true for $R=50$.
Rather, as $\alpha_1$ increases, $\uu_\gs$ tends to $\uu_\st$.
We would describe this behavior as shear thickening.

\begin{figure}
\centerline{(a)
\includegraphics[trim = 180pt 80pt 00pt 80pt,clip = true,width=0.45\textwidth]
{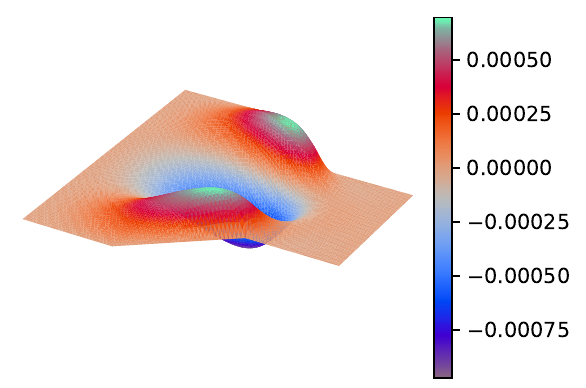}
\qquad     (b)
\includegraphics[trim = 180pt 80pt 00pt 80pt,clip = true,width=0.45\textwidth]
{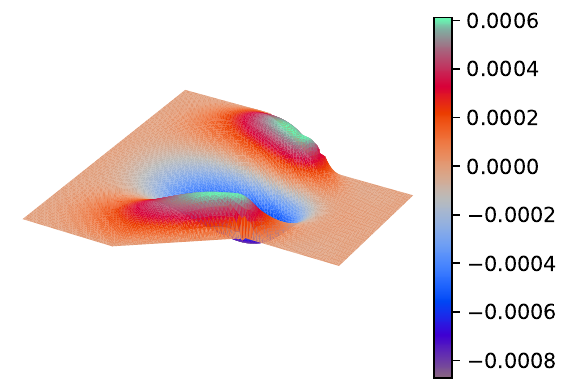}}
\caption{
Horizontal flow component of the difference $\uu_\gs-\uu_\st$ for
$\nu=1$, $U = 2^{-2}$ and (a) $\alpha_1=\alpha_2=0.02$, (b) $\alpha_1=\alpha_2=0.2$.
The computational domain $\domain$ is as specified in \eqref{eqn:expandom}, 
with $b_o = 1$, $b_i=1$, $L=1$, and  $H=0.5$.
The computational mesh was generated by four uniform refinements of the left-most 
mesh of figure \ref{fig:cduct-mesh}.}
\label{fig:gengtcontract}
\end{figure}

\begin{figure}
\centerline{(a)
\includegraphics[trim = 180pt 80pt 00pt 80pt,clip = true,width=0.45\textwidth]
{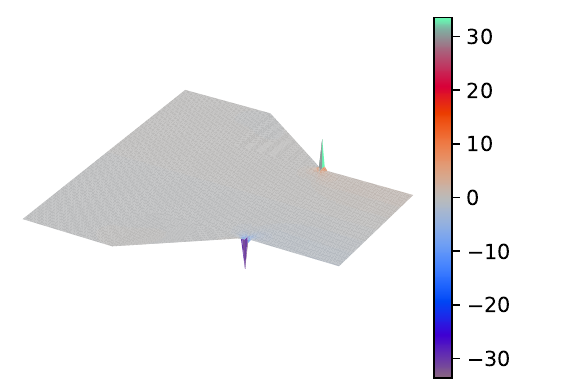}
\qquad     (b)
\includegraphics[trim = 180pt 80pt 00pt 80pt,clip = true,width=0.45\textwidth]
{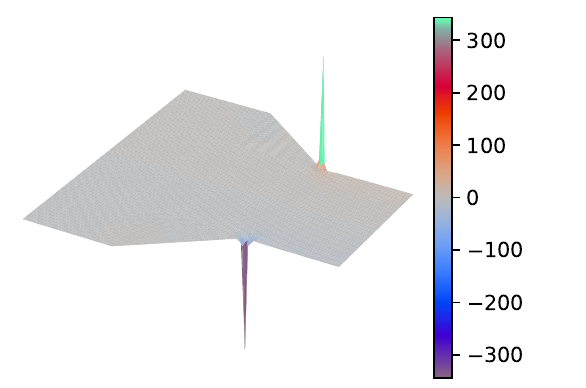}}
\caption{
Vertical component of the vector-valued auxiliary variable $\ww$ for
$\nu=1$, $U = 2^{-2}$ and (a) $\alpha_1=\alpha_2=0.02$, (b) $\alpha_1=\alpha_2=0.2$.
The computational domain $\domain$ is as specified in \eqref{eqn:expandom}, 
with $b_o = 1$, $b_i=1$, $L=1$, and  $H=0.5$.
The computational mesh was generated by four uniform refinements of the left-most 
mesh of figure \ref{fig:cduct-mesh}.}
\label{fig:wcontract}
\end{figure}

\subsection{Grade-two with $\alpha_1+\alpha_2\neq 0$}

In figure \ref{fig:gengtcontract} we contrast the difference in results between the
Stokes solution and the grade-two model when 
parameters $\alpha_1$ and $\alpha_2$ are are chosen independently, 
with $\alpha_1+\alpha_2\neq 0$.
We observe that we are no longer able to obtain solutions for
such large parameter values, and more care is required with defining the mesh.
In the plot on the right, with $\alpha_1 = \alpha_2 = 0.2$, we see some nonsmoothness 
in $\uu$ arising in the vicinity of the reentrant corners. The plot on the left 
with $\alpha_1 = \alpha_2 = 0.02$ remains smooth.
While we have $\norm{\uu^n}_{H^2(\Omega)}$ bounded in terms of 
$\norm{\ww^{n-1}}_{L^2(\Omega)}$, from \eqref{eqn:aldiventoo}-\eqref{eqn:diffrfsgeetoo}
the norm of the auxiliary variable $\norm{\ww^{n}}_{L^2(\Omega)}$ is bounded by
$\norm{\sdiv N(\uu^n,\pi^n)}_{L^2(\Omega)}$, which requires
$\norm{\uu^n}_{H^2(\Omega)}$ to be bounded.
More details on these bounds may be found in \cite{lrsBIBjd}.
Thus $\ww^n$ is sensitive to higher derivatives of $\uu^n$, and as illustrated in 
figure \ref{fig:wcontract}, these are not bounded for domains with nonconvex corners.
These corner singularities motiviate the localized refinement of the mesh for the 
accurate computation of the integral over the contraction boundary described in the
next section.  The localized refinement is further discussed in section 
\ref{sec:compudet}.

\section{Contraction rheometer}\label{sec:rheometer}

Contraction rheometers have been constructed e.g.~by Stading \cite{nystrom2017hyperbolic}.
The contraction zone generates a complex flow pattern that can be used to measure 
nonlinear relationships between the stress and rate of strain.
A recent paper \cite{lrsBIBix} examines the concept of identifiability 
for a rheometer for a given fluid model.
The typical experiment with a rheometer involves varying the flow rate and
measuring a force as a function of flow rate.
In \cite{lrsBIBix}, it is shown how to extract model parameters from such
a function for certain models.

For some models and rheometers, it is not possible to distinguish certain
parameters in a model, so these models are not identifiable by that rheometer.
For example \cite{lrsBIBix}, an extensional flow rheometer can determine the
sum $\alpha_1+\alpha_2$ for the grade-two model,
but it  does not distinguish the individual values $\alpha_i$.
And a simple shear rheometer is insensitive completely to both parameters $\alpha_i$.
Thus a natural question arises: can a contraction rheometer identify the
grade-two model?

\subsection{Force measurements}

The force $F$ that is measured by a contraction rheometer  is defined as: 
\begin{equation}\label{eqn:preforcint}
F = \int_{\bo} \nu \psi \, \nn^t (\nabla\uu + \nabla\uu^t) \hat{\xx} \, ds 
  - \int_{\bo} \psi \, p \, \nn^t\hat{\xx} \, ds, 
\end{equation}
where $\hat{\xx}$ is the first Euclidean basis vector and 
$\psi = 1$ for $x\in (0,L)$ and zero elsewhere (that is before and after the contraction,
see \eqref{eqn:expandom} 
and preceeding definitions).
Recall that $p$ is given by \eqref{eqn:peetopigoo}. 
Since $\uu$ is zero on the boundary, the formula simplifies to $p=\nu\pi$
on the support of $\psi$ on $\partial\Omega$.
Thus \eqref{eqn:preforcint} simplifies to
\begin{equation}\label{eqn:forcint}
F = \nu \bigg(\int_{\bo} \nu \psi \, \nn^t (\nabla\uu + \nabla\uu^t) \hat{\xx} \, ds 
  - \int_{\bo} \psi \, \pi \, \nn^t\hat{\xx} \, ds\bigg).
\end{equation}

Experimentally, this force is measured by a null balance device that
keeps the contraction portion of the device from moving.
The required force is thus proportional to $F$.
Computing $F(U)$ for various flow rates allows us to attempt to identify
the parameters of a model.
Typically, the flow rate $U$ is increased steadily from zero, possibly
plateauing at given values of $U$ temporarily to allow a steady state
to re-establish.
The limiting value of $F(U)/U$ for small values of $U$ is typically 
proportional to the viscosity $\nu$.
Other features of $F$ can be used to identify other parameters.

\omitit{
\subsection{Identifying identifiability}\label{subsec:ii}

Let us think abstractly about a general rheological model with parameters $c$.
What makes a model identifiable by a rheometer is that $F(U;c)$ is different
in some way for different values of the parameters $c$ of the model.
This can be determined computationally by computing $F(U;c)$ for different
values of $c$ to see if we get something different.

We may not get something different.
For example \cite{lrsBIBix}, an extensional flow rheometer gives the same values
$F(U)$ for all parameter values where the sum $\alpha_1+\alpha_2$ is fixed.
And a simple shear rheometer gives the same values for $F(U)$ for all 
values of $\alpha_i$.
Table \ref{tab:a1plusa2} shows results for the computed force 
integral $F(U)$ for varying values of $\alpha_1, \alpha_2$ and 
$\alpha_1 + \alpha_2$ and small values of $U$.
Here we may observe changing values of the computed force as $\alpha_1,
\alpha_2$ and $\alpha_1 + \alpha_2$ increasingly vary as the flow rate $U$ is increased. 
\begin{table}
\begin{center} 
\begin{tabular}{|l|c|c|c||c|c|}
\hline
$(\alpha_1,\alpha_2)$ & 
   $(0.3, 0.1)$& 
   $(0.2, 0.2)$&  
   $(0.1, 0.3)$& 
   $(0.2, 0.1)$& 
   $(0.1, 0.2)$\\
\hline
$U = 2^{-4} $ &-11.35767 &-11.32156 &-11.28623 &-11.25167 &-11.21633 \\
$U = 2^{-5} $ &-11.15493 &-11.13734 &-11.11977 &-11.10249 &-11.08495 \\
$U = 2^{-6} $ &-11.05447 &-11.04573 &-11.03700 &-11.02835 &-11.01962 \\
\hline
\end{tabular}
\end{center}
\vspace{-5mm}
\caption{$F(U)/U$ 
with force $F$ defined in \eqref{eqn:forcint}, over the computational domain
$\domain$ is as specified in \eqref{eqn:expandom}, with $\nu=1$,
$b_i=1$, $b_o=1$, $L=1$, $H=0.5$.
Results for varying values of $U$, $\alpha_1$ and
$\alpha_2$ show the force integral $F$ varies with $\alpha_1$ and $\alpha_2$, even
where their sum remains constant.
Computations were performed on a local refinement of the mesh show in 
figure \ref{fig:cduct-mesh}, as descripted in subsection \ref{subsec:locref}
} 
\label{tab:a1plusa2}
\end{table}

Thus a simple approach to the identifiability of the grade-two model by
a contraction rheometer is to plot $F(U)$, or $F(U)/U$, for some interval
of $U$ values, for various values of $\nu$ and $\alpha_i$.
If we plot $F(U)/U$, then it is plausible that the values for small $U$ will be
independent of $\alpha_i$ and proportional to $\nu$,
as this should correspond to Stokes flow,
and this expected result gives an internal check.
Results for this internal check for the normalized force integral $F(U)/U$ are shown in 
table \ref{tab:aitkendenti},
where we see for small values of $U$ the force integral $F(U)$ satisfies
the predicted relation.
The computations for tables 
\ref{tab:a1plusa2} and \ref{tab:aitkendenti}
 were performed on 
a locally refined mesh with the refinements focused on the contraction boundary
and its endpoints. Specific details on the computation and the mesh are given in section 
\ref{sec:compudet} and subsection \ref{subsec:locref}.
}

\begin{table}
\begin{center} 
\begin{tabular}{|c|c|c|c||c|c|c|}
\hline
viscosity
&  \multicolumn{3}{c||}{$\nu = 1$}
&  \multicolumn{3}{c|}{$\nu = 2$}\\
\hline
   $(\alpha_1,\alpha_2)$ & 
   $(0.1, 0.1)$& $(0.1, 0.2)$&$(0.2, 0.1)$& 
   $(0.1, 0.1)$& $(0.1, 0.2)$&$(0.2, 0.1)$\\
\hline
$U = 2^{-6}$ &-11.00226 &-11.01962 &-11.02835 &-21.95676 &-21.97412 &-21.98282 \\
$U = 2^{-7}$ &-10.97838 &-10.98706 &-10.99141 &-21.93292 &-21.94159 &-21.94593 \\
$U = 2^{-8}$ &-10.96646 &-10.97080 &-10.97297 &-21.92101 &-21.92534 &-21.92751 \\
\hline
$\Delta^2  $ &-10.95458 &-10.95458 &-10.95458 &-21.90911 &-21.90912 &-21.90913 \\
\hline
\end{tabular}
\end{center}
\vspace{-5mm}
\caption{$F(U)/U$ with force $F$ defined in \eqref{eqn:forcint}, 
over the computational domain $\domain$ as specified in \eqref{eqn:expandom}, 
with $b_i=1$, $b_o=1$, $L=1$, $H=0.5$.
Results are shown for small values of $U$ and varying values of $\alpha_1$ and 
$\alpha_2$ with $\nu = 1$ and $\nu = 2$. The bottom row shows the extrapolated
values as $U \rightarrow 0$ computed by the Aitken $\Delta^2$ process using the three
rows above.  
Computations were performed on a local refinement of the mesh show in 
figure \ref{fig:cduct-mesh}, as descripted in subsection \ref{subsec:locref}
} 
\label{tab:aitkendenti}
\end{table}

In particular, the bottom row of table \ref{tab:aitkendenti}, which was computed by
applying Aitken's $\Delta^2$ process to the two rows above to compute a value of 
$F(U)/U$ as $U \rightarrow 0$, shows the computed force is approximately $-10.9546\nu$.
We next turn our attention to the identifiability of $\alpha_1$ and $\alpha_2$.

\subsection{Visualizing the data}\label{subsec:vis}
In this section we consider the computed data from our force measurement
simulations 
$ f(U,\nu,\balpha)=F(U)/U$ for $\balpha=(\alpha_1,\alpha_2)$.
Writing $\balpha$ in polar coordinates as 
$\balpha = \alpha(\cos(\theta),\sin(\theta))$, we can visualize the dependence of
$f$ on $U$, $\alpha$ and $\theta$ as follows. This allows us to better understand
regimes where the parameters $\alpha_1, \alpha_2$ could be identifiable using
this rheometer.

First in figure \ref{fig:dataU-thetadep} we consider snapshots in $U$ while $\theta$ is 
varied for different values of $\alpha$.
\begin{figure}
\centering
\includegraphics[trim = 0pt 0pt 0pt 0pt,clip = true,width=0.32\textwidth]
{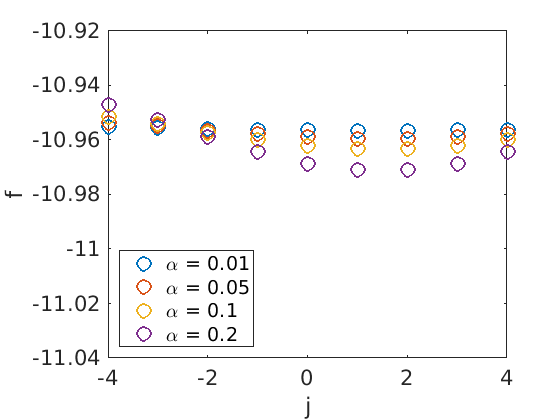}
\includegraphics[trim = 0pt 0pt 0pt 0pt,clip = true,width=0.32\textwidth]
{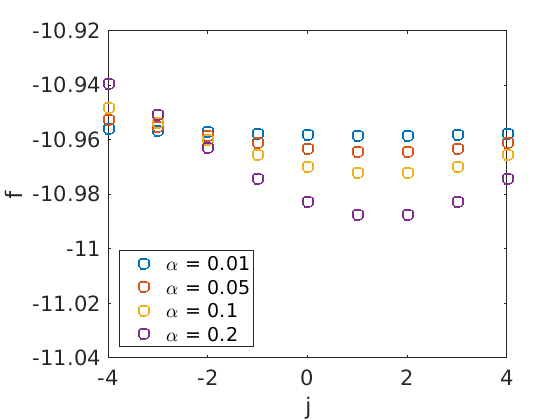}
\includegraphics[trim = 0pt 0pt 0pt 0pt,clip = true,width=0.32\textwidth]
{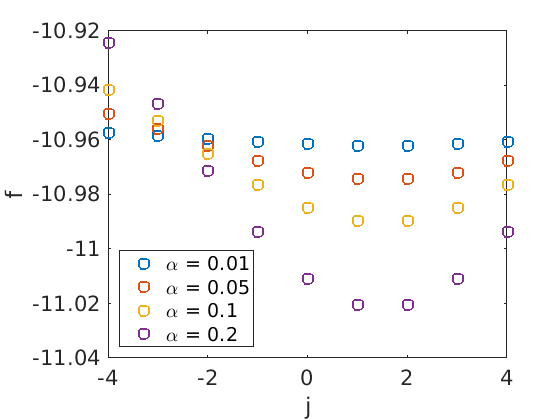}
\caption{Generated data $f(U,\nu,\balpha) = F(U)/U$ with the same computational domain
and mesh as in table \ref{tab:aitkendenti} for
$\balpha = \alpha (\cos( j \pi/8),\sin(j \pi/8))$ with $\nu = 1$ and varying $\alpha$
and $j$.
Left: $U = 2^{-8}$; center: $U = 2^{-7}$; right: $U = 2^{-6}$.
}
\label{fig:dataU-thetadep}
\end{figure}
We make three important observations from these three snapshots.  First, $f$ has a 
nonlinear dependence on $\theta$ for fixed $U$. 
In contrast, as shown in 
figure \ref{fig:polyerrU} and table \ref{tab:maxIerr} 
which display the difference between the computed data and the best fit linear, quadratic
and cubic polynomials for a range of $U$ values with $\alpha$ and $\theta$ fixed,
$f$ is close to linear in $U$;
this is further illustrated in 
figures \ref{fig:data-theta-cross} and \ref{fig:data-reflect}.
Second, it appears that each of the 
fixed-$\alpha$ trajectories  of figure \ref{fig:dataU-thetadep}
cross at $\theta = j \pi/8$ for $j \approx -2.5$,
which we will look at more closely in table \ref{tab:alcross} and see further illustrated
in figure \ref{fig:data-theta-cross}.  
Third, it appears $f$ is symmetric 
over $\theta \approx j \pi/8$ for $j = 1.5$, which we will look at more closely in figure
\ref{fig:data-reflect}. 

To visualize the near linearity of  $f$ as a function of $U$, we computed the 
best fit regression line $\cI_1 f$, quadratic fit $\cI_2 f$, and cubic fit $\cI_3f$ to the
10 data points for $U = 0.01, 0.02, \ldots, 0.1,$ with fixed 
$\theta \in \{\pm 3 \pi/16, \pm \pi/16\}$ and $\alpha \in \{0.03, 0.07, 0.14\}$.
The results of the errors $f - \cI f$ normalized by its maximum value for
each line are shown in figure \ref{fig:polyerrU}.
The normalization factors are shown in table \ref{tab:maxIerr}.
We see that
\begin{itemize}
\item
the error in each linear regression
line is approximately quadratic, with its maximum magnitude on the order of $10^{-6}$; 
\item
the error in each best-fit quadratic is approximately cubic, 
with its maximum magnitude on the order of $10^{-8}$;
\item
 and the error in the best-fit cubic is not clearly structured but 
has maximum magnitude on the order of $10^{-10}$. 
\end{itemize}

\begin{figure}
\centering
\includegraphics[trim = 0pt 0pt 0pt 0pt,clip = true,width=0.30\textwidth]
{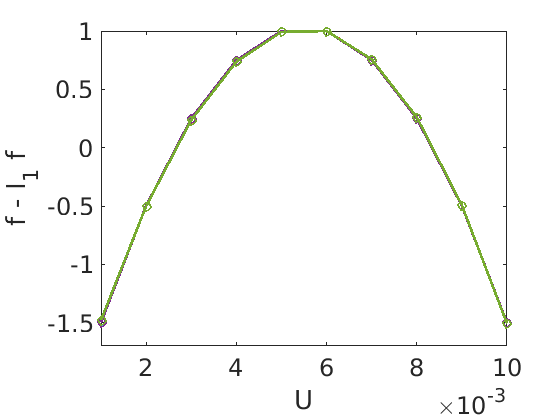}
~\hfil~
\includegraphics[trim = 0pt 0pt 0pt 0pt,clip = true,width=0.30\textwidth]
{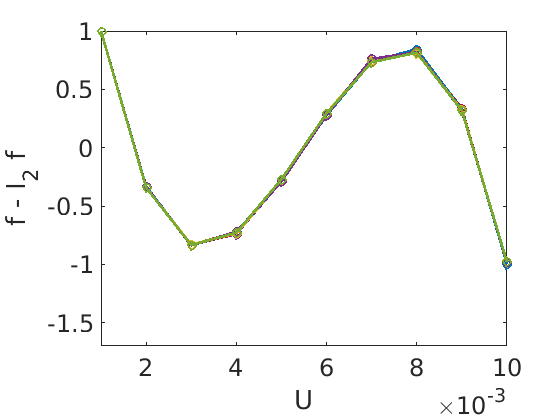}
~\hfil~
\includegraphics[trim = 0pt 0pt 0pt 0pt,clip = true,width=0.30\textwidth]
{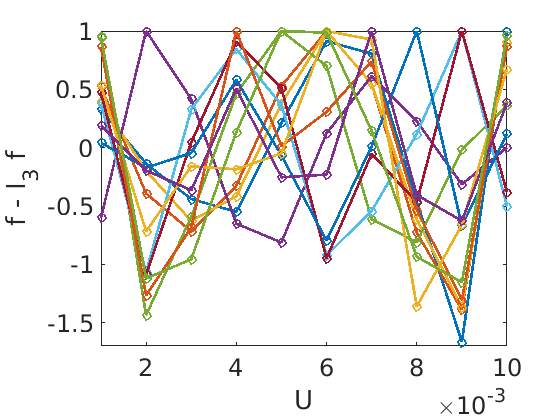}
\caption{Normalized difference of the error in the least squares regression line (left), quadratic fit (center) and cubic fit (right) and the generated data 
$f(U,\nu,\balpha) = F(U)/U$. The normalization factors (maximum errors) are shown in 
table \ref{tab:maxIerr}. Each plot shows an overlay of 12 $\alpha,\theta$ pairs for
$\balpha = \alpha (\cos(\theta),\sin(\theta))$ with 
$\theta \in\{ \pm \pi/16, \pm 3 \pi/16\}$, $\alpha \in\{ 0.03, 0.07, 0.14\}$ 
and $U$ increasing from $0.01$ to $0.1$.
The data $f(U,\nu,\balpha)$ were generated with the same 
computational domain
and mesh as in table \ref{tab:aitkendenti}. 
}
\label{fig:polyerrU}
\end{figure}

\begin{table}
\begin{tabular}{|l||c | c | c | c |}
\hline
& \multicolumn{4}{c|}{max $(f - \cI_1 f)$}  \\ \hline
\hline
$\alpha / \theta$ &$-3 \pi/16$ & $-\pi/16$ & $\pi/16$ & $3 \pi/16$\\
\hline
$0.03$ & 3.810e-6 & 0.419e-6 & 0.533e-6 & 0.563e-6 \\
$0.07$ & 0.576e-6 & 1.052e-6 & 1.360e-6 & 1.420e-6 \\
$0.14$ & 1.306e-6 & 2.532e-6 & 3.339e-6 & 3.423e-6\\
\hline
\end{tabular}
\begin{tabular}{|l||c | c | c | c |}
\hline
\multicolumn{4}{|c|}{max $(f - \cI_2 f)$}  \\ \hline
\hline
$-3 \pi/16$ & $-\pi/16$ & $\pi/16$ & $3 \pi/16$\\
\hline
5.760e-8 & 0.020e-8 & 0.040e-8 & 0.043e-8  \\
0.054e-8 &0.305e-8 & 0.562e-8 & 0.579e-8  \\
0.554e-8 &2.622e-8 & 4.677e-8 & 4.769e-8  \\
\hline
\end{tabular} \\
\begin{tabular}{|l||c | c | c | c |}
\hline
& \multicolumn{4}{c|}{max $(f - \cI_3 f)$}  \\ \hline
\hline
$\alpha / \theta$ &$-3 \pi/16$ & $-\pi/16$ & $\pi/16$ & $3 \pi/16$\\
\hline
$0.03$ &6.385e-10 &0.052e-10 &0.049e-10 & 0.038e-10 \\
$0.07$ &0.062e-10 &0.057e-10 &0.057e-10 & 0.097e-10 \\
$0.14$ &0.124e-10 &0.694e-10 &1.283e-10 & 0.564e-10 \\
\hline
\end{tabular}
\caption{Maximum error between the computed data points $f(U,\nu,\balpha) = F(U)/U$
and the regression line $\cI_1 f$, quadratic fit $\cI_2 f$, and cubic fit $\cI_3 f$. Each line fits the data for
the flow rate $U$ increasing linearly from $0.01$ to $0.1$ with fixed values of 
$\alpha = \{0.03, 0.07, 0.14\}$ and $\theta = \{\pm \pi/16, \pm 3 \pi/16\}$.
}
\label{tab:maxIerr}
\end{table}

\begin{table}
\begin{tabular}{|l||c|c||c|c||c|c||}
\hline 
   & \multicolumn{2}{c||}{$U = 2^{-8},\theta_8 = -2.50266827$} 
   & \multicolumn{2}{c||}{$U = 2^{-7},\theta_7 = -2.50268795$} 
   & \multicolumn{2}{c||}{$U = 2^{-6},\theta_6 = -2.50272560$} \\
\hline
$\alpha$ & $\theta= \theta_8 \cdot \pi/8$ & $\theta = -2.5 \pi/8$ 
         & $\theta= \theta_7 \cdot \pi/8$ & $\theta = -2.5 \pi/8$ 
         & $\theta= \theta_6 \cdot \pi/8$ & $\theta = -2.5 \pi/8$ \\
\hline
0.01     & -10.9556149 & -10.9556158 & -10.9566795 & -10.9566812
         & -10.9588091 & -10.9588125\\
0.05     & -10.9556149 & -10.9556190 & -10.9566795 & -10.9566878
         & -10.9588091 & -10.9588259 \\
0.1      & -10.9556149 & -10.9556231 & -10.9566795 & -10.9566961
         & -10.9588091 & -10.9588427\\
0.2      & -10.9556149 & -10.9556313 & -10.9566795 & -10.9567126
         & -10.9588091 & -10.9588762\\
\hline
\end{tabular}
\caption{Generated data $f(U,\nu,\balpha) = F(U)/U$ with the same computational domain
and mesh as in table \ref{tab:aitkendenti}.
Results are shown for values of $\balpha = \alpha (\cos(\theta),\sin(\theta))$
near the crossing of trajctories of constant $\alpha$ values shown in figure
\ref{fig:dataU-thetadep}. The crossing for each value of $U$ occurs at 
approximately $\theta = \theta_k \cdot \pi/8$, $k = 6,7,8$.
Left: $U = 2^{-8}$; center: $U = 2^{-7}$; right: $U = 2^{-6}$.
} 
\label{tab:alcross}
\end{table}
\begin{figure}
\centering
\includegraphics[trim = 0pt 0pt 0pt 0pt,clip = true,width=0.320\textwidth]
{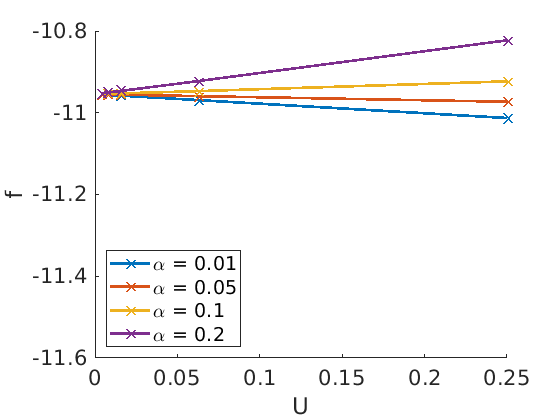}
~\hfil~
\includegraphics[trim = 0pt 0pt 0pt 0pt,clip = true,width=0.30\textwidth]
{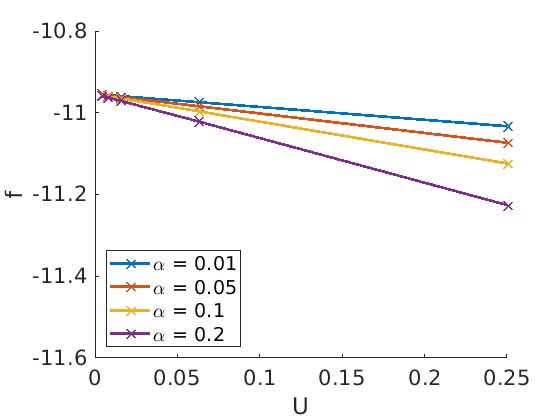}
~\hfil~
\includegraphics[trim = 0pt 0pt 0pt 0pt,clip = true,width=0.30\textwidth]
{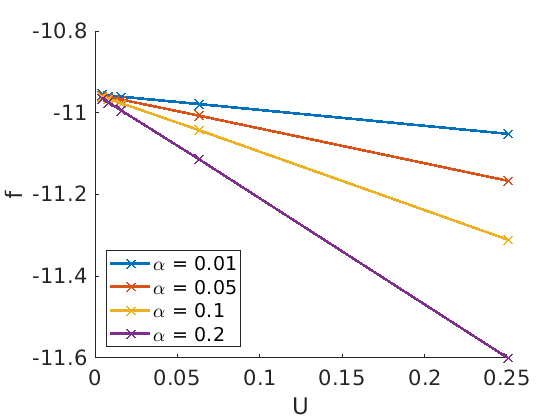}
\caption{Generated data $f(U,\nu,\balpha) = F(U)/U$ with the same computational domain
and mesh as in table \ref{tab:aitkendenti} showing how $f$ changes with $U$ for
$U = 2^{-j}, j = \{8,7,6,4,2\}$ and different values of $\alpha$, the magnitude of 
$\balpha = \alpha(\cos(\theta),\sin(\theta))$.
Left: $\theta = -3 \pi/8$; center: $\theta = - 2 \pi /8$; right: $\theta = - \pi/8$.
}
\label{fig:data-theta-cross}
\end{figure}

To investigate the crossing of trajectories of each fixed value of $\alpha$ near
$\theta = -5 \pi/16$, for each $U_k = 2^{-k}$, $k = 6,7,8$, we computed the 
intersection of the linear interpolants connecting
$f$ at $\theta_L = -2.65 \pi/8$ (left of the crossing) to 
$\theta_R = -2.5 \pi/8$ (right of the crossing) for $\alpha = 0.01$ and $\alpha = 0.2$.
The computed crossing points $ \theta = \theta_k \pi/8$ are shown in table 
\ref{tab:alcross}
along with the value of $f(U,1,\balpha)$ at each of the given coordinates.
As shown in the table, $f$ is constant through 7 decimal
places at each $\theta = \theta_k \pi/8, U= 2^{-k}$ pair of coordinates (not shown: 
the values are decreasing with $\alpha$ for each at the 9th digit).  
At $\theta = -2.5 \pi/8$, about
$0.001$ to the right in $\theta$, the values of $f$ are increasing with $\alpha$, 
where the differences are seen in the fifth decimal place, 
which is on the order of the tolerance of our solve.
This shows that $\alpha$, the magnitude of $\balpha$, is not identifiable using 
the contraction rheometer if the argument $\theta$ of $\balpha$ is close to
$-2.5 \pi/8$.

Figure \ref{fig:data-theta-cross} shows snapshots of $f$ as a function of $U$
for different values of $\alpha$ with $\theta = -3\pi/8$ (left) and $\theta = -2\pi/8$
(center), which lie on either side of the crossing.  We observe that the
constant-$\alpha$ trajectories are in opposite order on these two plots.  
Comparing the center plot for $\theta = -2 \pi/8$ and the right plot for 
$\theta = -\pi/8$, we can see how the angle between fixed-$\alpha$ trajectories changes
for different values of $\theta$.

\begin{figure}
\centering
\includegraphics[trim = 0pt 0pt 0pt 0pt,clip = true,width=0.30\textwidth]
{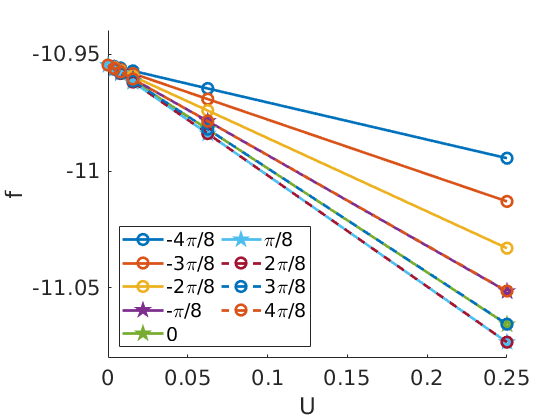}
~\hfil~
\includegraphics[trim = 0pt 0pt 0pt 0pt,clip = true,width=0.30\textwidth]
{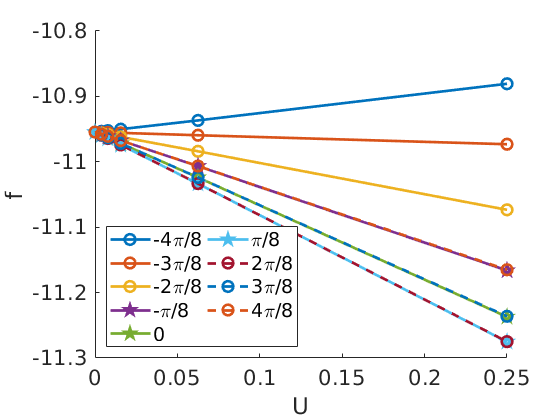}
~\hfil~
\includegraphics[trim = 0pt 0pt 0pt 0pt,clip = true,width=0.30\textwidth]
{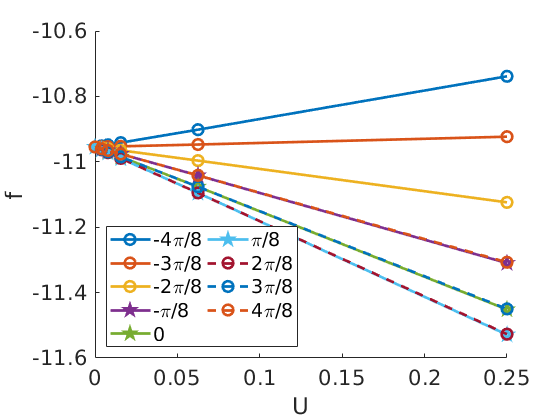}
\caption{Generated data $f(U,\nu,\balpha) = F(U)/U$ with the same computational domain
and mesh as in table \ref{tab:aitkendenti} showing how $f$ changes with $U$ for
$U = 2^{-j}, j = \{8,7, 6,4,2\}$ and different values of $\theta$, the argument of 
$\balpha = \alpha(\cos(\theta),\sin(\theta))$.
Left: $\alpha = 0.01$; center: $\alpha = 0.05$; 
right: $\alpha = 0.1$.
}
\label{fig:data-reflect}
\end{figure}
\begin{figure}
\centering
\includegraphics[trim = 0pt 0pt 0pt 0pt,clip = true,width=0.30\textwidth]
{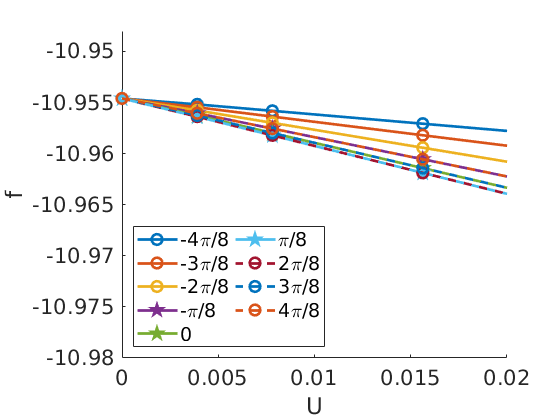}
~\hfil~
\includegraphics[trim = 0pt 0pt 0pt 0pt,clip = true,width=0.30\textwidth]
{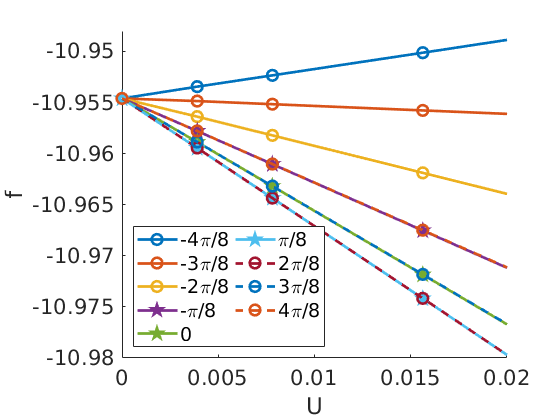}
~\hfil~
\includegraphics[trim = 0pt 0pt 0pt 0pt,clip = true,width=0.30\textwidth]
{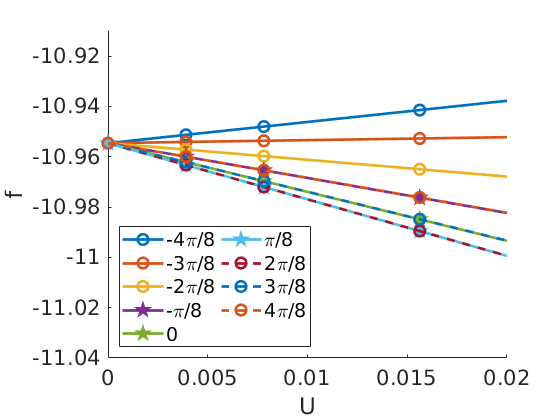}
\caption{Zoomed in view of figure \ref{fig:data-reflect} showing generated data $f(U,\nu,\balpha) = F(U)/U$ with the same computational domain
and mesh as in table \ref{tab:aitkendenti} showing how $f$ changes with $U$ for
$U = 2^{-j}, j = \{8,7,6,4,2\}$ and different values of $\theta$, the argument of 
$\balpha = \alpha(\cos(\theta),\sin(\theta))$.
Left: $\alpha = 0.01$; center: $\alpha = 0.05$; 
right: $\alpha = 0.1$. 
The first two plots use the same scaling.
}
\label{fig:data-reflect-zoom}
\end{figure}

Our third observation from figure \ref{fig:dataU-thetadep} is that $f$ appears 
symmetric around $\theta_S \approx 3 \pi/16$, that is, 
\begin{equation}\label{eqn:symmtheta}
f(U,\nu,\alpha(\cos(\theta_S+t),\sin(\theta_S+t))
=f(U,\nu,\alpha(\cos(\theta_S-t),\sin(\theta_S-t)),
\end{equation}
where the equality holds up to the tolerance of our solve.
This matching of trajectories of $f$ as $U$ is increased is further illustrated
in figure \ref{fig:data-reflect} which shows $f$
for four different values of $\alpha$ and varying $U$.
The fixed-$\theta$ trajectories are nearly linear in $U$; and, 
\begin{itemize}
\item the data for $\theta = \pi/2=8\pi/16$ overlays the data for 
$\theta = -\pi/8=-2\pi/16$ ($t=5\pi/16$ in \eqref{eqn:symmtheta}),

\item the data for $\theta = 3\pi/8=6\pi/16$ overlays the data for $\theta = 0$
($t=3\pi/16$ in \eqref{eqn:symmtheta}), and

\item the data for $\theta = \pi/4=2\pi/8=4\pi/16$ overlays the data
for $\theta = \pi/8=2\pi/16$
($t=\pi/16$ in \eqref{eqn:symmtheta}).
\end{itemize}
Figure \ref{fig:data-reflect-zoom} shows a zoomed-in view of figure 
\ref{fig:data-reflect}. 
The point at $(0,0)$ was added to each of these plots to illustrate again 
the near-linearity in $U$ of data $f$ for small flow rates $U$.
In terms of identifying $\balpha$ given a sequence of measurements of $f(U,\nu,\balpha)$ 
for varying flow rates $U$, we conclude that we can only expect to identify  
$\balpha = \alpha (\cos(\theta),\sin(\theta))$ in a limited range in arguments 
$\theta$ of $\balpha$, namely $\theta \in (-5 \pi/16, 3 \pi/16]$.

\section{Computational details}\label{sec:compudet}

For our simulations of flow in a contracting duct, as described in section 
\ref{sec:cduct-flow}, we solved to a default tolerance of 
$\nr{\uu^n - \uu^{n-1}}< 10^{-5}$, and an
iterated-penalty-method (IPM) \cite{lrsBIBgd,lrsBIBih}
tolerance of 
$\nr{\sdiv (\uu^{n,l}) }< 10^{-10}$, with IPM parameter $\rho = 10^4$.
We used Lagrange degree $k=4$ elements for the space $V_h$ from subsection 
\ref{subsec:disc-alg}. 
All computations were performed using FEniCS 
\cite{FEniCSbook,lrsBIBih}.

\subsection{Localized mesh refinement}\label{subsec:locref}
For the contraction rheometer computations described in Section \ref{sec:rheometer},
we performed a local refinement of the left-hand mesh in figure 
\ref{fig:cduct-mesh}. The localized refinements consisted of boundary refinements 
(denoted $r_b$) and point refinements (denoted $r_p$). 
The boundary refinements first mark any element sharing an edge with the 
contraction boundary. These elements were identified by (a) having a maximum distance of 
$0.65$ to the points $(0.4, \pm 0.5)$, and (b) sharing an edge with the boundary.
The refinement includes subdivision of the marked elements, followed by a  completion 
step refining certain neighboring elements to ensure the mesh is conforming. 
The point refinements (denoted $r_p$) first
subdivide any element that contains one of the four endpoints of the contraction 
boundary, namely $(1, \pm 0.5)$ and $(0,\pm 1)$; 
followed by a completion step to ensure the mesh is conforming.
Meshes with $r_b = 2, r_p = 0$ and $r_b=2, r_p = 4$ are illustrated in the
center and right of figure \ref{fig:cduct-mesh}.

\begin{figure}
\includegraphics[trim = 40pt 40pt 40pt 40pt,clip = true,width=0.30\textwidth]
{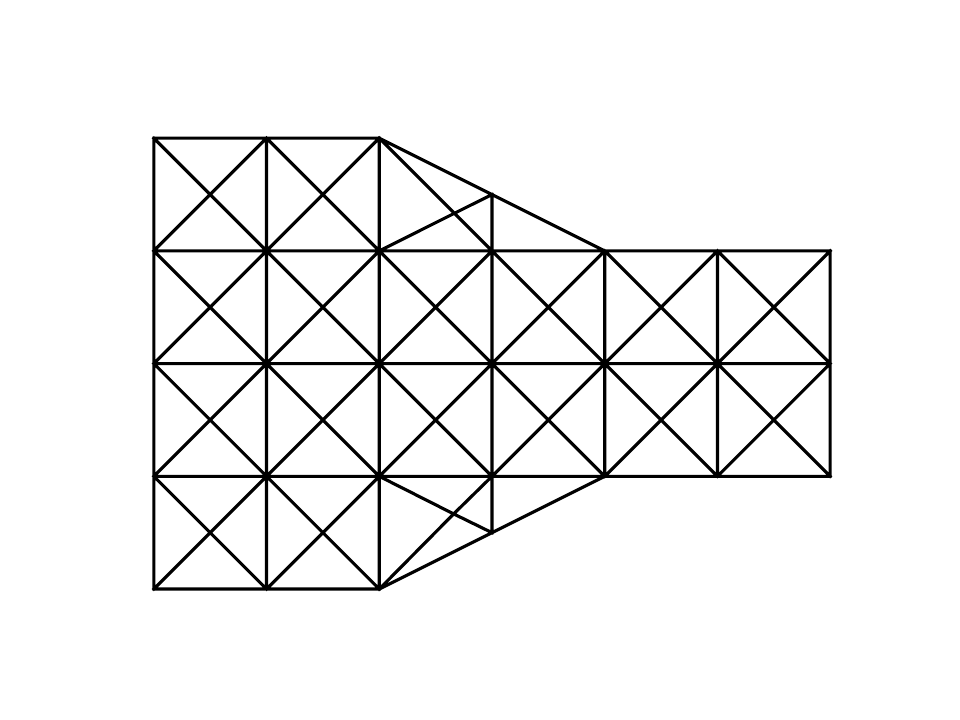}
\includegraphics[trim = 40pt 40pt 40pt 40pt,clip = true,width=0.30\textwidth]
{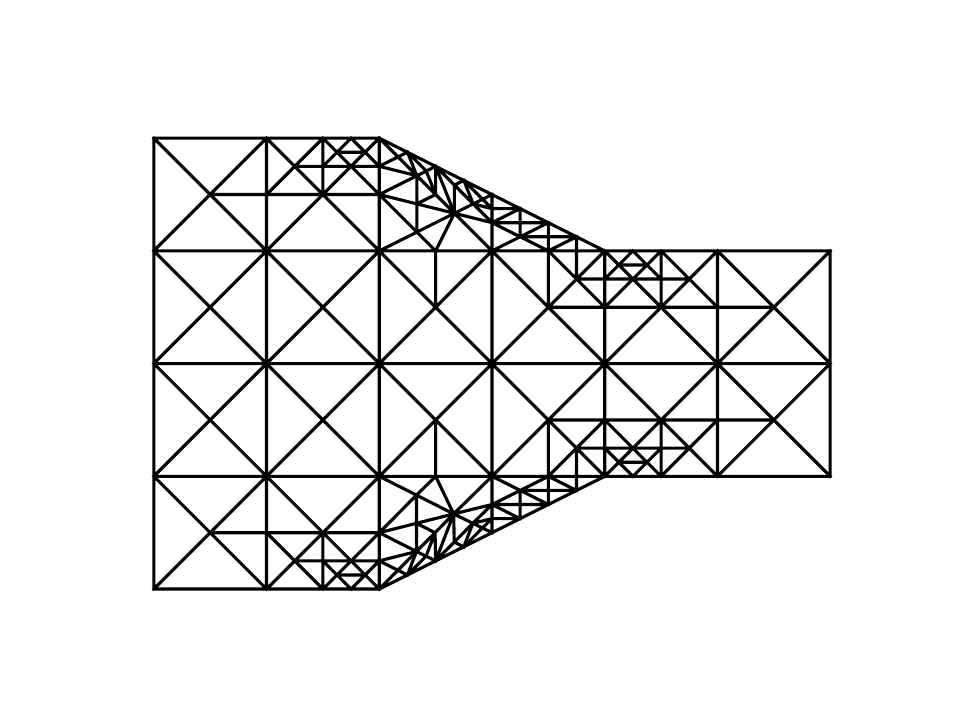}
\includegraphics[trim = 40pt 40pt 40pt 40pt,clip = true,width=0.30\textwidth]
{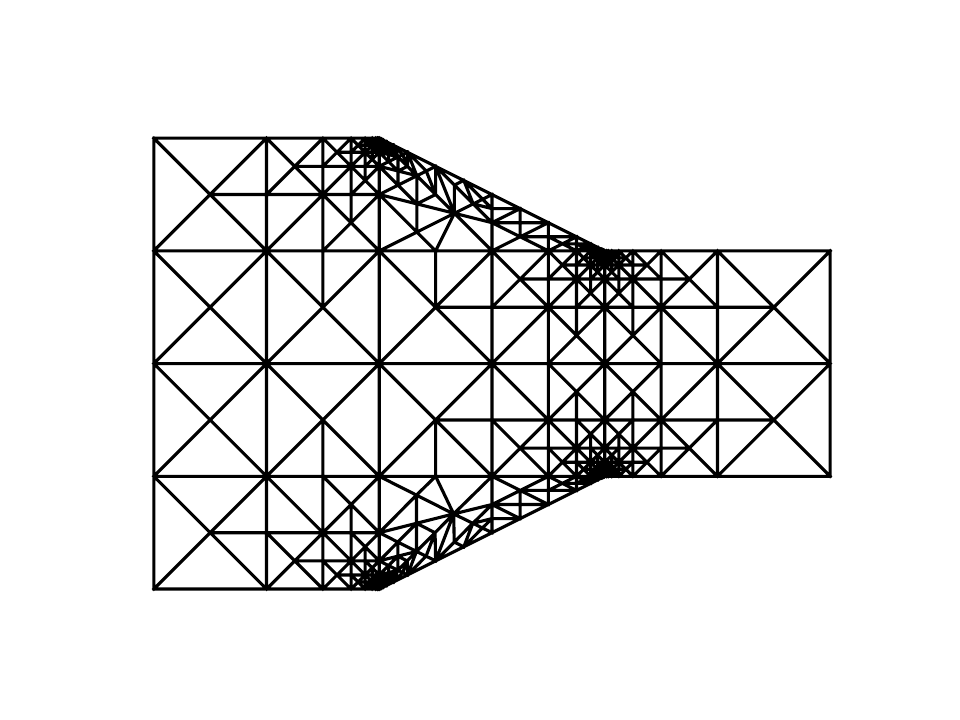}
\caption{Boundary refinements ($r_b$) and point refinements ($r_p$) of 
the computational mesh for the contracting duct.
Left: $r_b = 0, r_p = 0$; center: $r_b = 2, r_p = 0$; right: $r_b = 2, r_p = 4$.}
\label{fig:cduct-mesh}
\end{figure}

To determine a computationally efficient but accurate mesh, we first compared 
uniform refinements (denoted $r_u$) to refinements focused on the contraction boundary,
as described above.  As shown in table \ref{tab:mesh-conv0p}, without any point 
refinements, boundary refinements 
consistently gave the same result as uniform refinements to two decimal places in the 
computation of the normalized force integral $F/U$, while neither converged to even a 
single decimal place after seven refinements. 

As shown in table \ref{tab:mesh-conv12p}, with 12 point refinements, uniform and boundary
refinement gave the same results to three decimal places on the least refined and four
decimal places on the most refined meshes. By comparison with table 
\ref{tab:mesh-conv4pt}, we see that after 7 refinements of
either type, $F/U$ converged to two decimal places.
To maintain efficiency, we eliminated the uniform refinements, and attained 
significantly more accuracy using a combination of point refinements and boundary 
refinements, as shown in 
table \ref{tab:mesh-conv4pt}.  In particular, 9 boundary refinements followed by
12 point refinements resolves $F/U$ to at least $5 \times 10^{-4}$.  
Increasing $r_p$ from 12
to 16 however did not increase the accuracy of the computation. The results of subsection
\ref{subsec:vis} 
were computed with the most accurate
combination found here: $r_b = 9$ and $r_p = 12$, for a total of 315,922 degrees of 
freedom in the discrete space $V_h$, using quartic (degree 4) vector-valued Lagrange 
basis functions.
 
\begin{table}
\begin{center} 
\begin{tabular}{|l|r|r|r|r|r|r|r|r|}
\hline
$r_p = 0$ &$r_b = 0$ & $r_b = 1$ & $r_b = 2$ & $r_b = 3$& $r_b = 4$& $r_b = 5$ 
& $r_b = 6$ & $r_b = 7$\\ 
\hline
$r_u = 0$ & -0.00000 & -5.08710 & -7.80120 & -9.30473 &-10.15756 &-10.64941 &-10.93652 & -11.10571 \\
$r_u = 1$ & -5.08429 & -7.80074 & -9.30503 &-10.15792 &-10.64980 &-10.93692 &-11.10612 &           \\  
$r_u = 2$ & -7.80075 & -9.30300 &-10.15784 &-10.64975 &-10.93688 &-11.10607 &          & \\  
$r_u = 3$ & -9.30527 &-10.15532 &-10.64951 &-10.93667 &-11.10587 &          &          & \\
$r_u = 4$ &-10.15816 &-10.64720 &-10.93639 &-11.10560 &          &          &          & \\ 
\hline
\end{tabular}
\end{center}
\vspace{-5mm}
\caption{$F/U$ with force $F$ defined in \eqref{eqn:forcint}, 
over the computational domain
$\domain$ is as specified in \eqref{eqn:expandom}, with $b_i=1$, $b_o=1$, $L=1$, $H=0.5$
and parameters $\alpha_1= 0.3$, $\alpha_2=0.1$, $\nu = 1$ and $U = 2^{-4}$.
The results shown for $r_u$ uniform refinements and $r_b$ boundary refinements
of the mesh shown in figure \ref{fig:cduct-mesh} demonstrate
that without point refinements, uniform refinements can be exchanged for refinements 
over the contraction boundary to increase the accuracy of the computation with fewer 
degrees of freedom.
} 
\label{tab:mesh-conv0p}
\end{table}

\begin{table}
\begin{center} 
\begin{tabular}{|l|r|r|r|r|r|r|r|r|}
\hline
$r_p = 0$ &$r_b = 0$ & $r_b = 1$ & $r_b = 2$ & $r_b = 3$& $r_b = 4$& $r_b = 5$ 
& $r_b = 6$ & $r_b = 7$\\ 
\hline
$r_u = 0$ &-11.33672 &-11.34446  &-11.34942  &-11.35244 &-11.35434 &-11.35553 &-11.35626 &-11.35670 \\ 
$r_u = 1$ &-11.34479 &-11.34947  &-11.35249  &-11.35439 &-11.35558 &-11.35631 &-11.35674 &          \\
$r_u = 2$ &-11.34960 &-11.35246  &-11.35435  &-11.35554 &-11.35627 &-11.35671 &          &          \\ 
$r_u = 3$ &-11.35258 &-11.35434  &-11.35553  &-11.35626 &-11.35669 &          &          &          \\
$r_u = 4$ &-11.35440 &-11.35552  &-11.35625  &-11.35668 &          &          &          &          \\ 
\hline
\end{tabular}
\end{center}
\vspace{-5mm}
\caption{$F/U$ with force $F$ defined in \eqref{eqn:forcint}, 
over the computational domain
$\domain$ is as specified in \eqref{eqn:expandom}, with $b_i=1$, $b_o=1$, $L=1$, $H=0.5$
and parameters $\alpha_1= 0.3$, $\alpha_2=0.1$, $\nu = 1$ and $U = 2^{-4}$.
The results shown for $r_u$ uniform refinements, $r_b$ boundary refinements and 
12 point refinements
of the mesh shown in figure \ref{fig:cduct-mesh} demonstrate
that used together with point refinements, uniform refinements can be exchanged for refinements over the 
contraction boundary to increase the accuracy of the computation with fewer degrees of freedom.
} 
\label{tab:mesh-conv12p}
\end{table}

\begin{table}
\begin{center} 
\begin{tabular}{|l|r|r|r|r||r|}
\hline
$r_u = 0$ &$r_p = 0$ & $r_p = 4$& $r_p = 8$& $r_p = 12$ & $r_p = 16 $  \\ 
\hline
$r_b = 4$ &-11.35705 &-11.35434 &-11.33619 &-11.20465 &-10.15756 \\ 
$r_b = 5$ &-10.64941 &-11.26481 &-11.34434 &-11.35553 &-11.35706 \\
$r_b = 6$ &-10.93652 &-11.30100 &-11.34936 &-11.35626 &-11.35716 \\ 
$r_b = 7$ &-11.10571 &-11.32291 &-11.35247 &-11.35670 &-11.35709 \\
$r_b = 8$ &-11.20614 &-11.33624 &-11.35439 &-11.35738 &-11.35547 \\
$r_b = 9$ &-11.26616 &-11.34440 &-11.35559 &-11.35767 &-11.35522 \\
\hline
\end{tabular}
\hfil
\end{center}
\vspace{-5mm}
\caption{$F/U$ with force $F$ defined in \eqref{eqn:forcint}, 
over the computational domain
$\domain$ is as specified in \eqref{eqn:expandom}, with $b_i=1$, $b_o=1$, $L=1$, $H=0.5$
and parameters $\alpha_1= 0.3$, $\alpha_2=0.1$, $\nu = 1$ and $U = 2^{-4}$.
The results shown for no uniform refinements, $r_b$ boundary refinements and
$r_p$ point refinements
of the mesh shown in figure \ref{fig:cduct-mesh} demonstrate that sufficient refinement
at the four endpoints of the contraction boundary together with boundary refinements 
provides an increase in accuracy.  As seen with $r_p = 16$, too many refinements at
those points does not further increase the accuracy.
} 
\label{tab:mesh-conv4pt}
\end{table}

\subsection{Computational mesh}\label{subsec:computmesh}

In this section we summarize our findings on defining an appropriate
computational mesh for this problem. 
First, the discrete inf-sup constant should not degenerate as the mesh size is
decreased. 
The convergence of the IPM algorithm slows with diminishing returns if the inf-sup
constant is sufficiently small. This slow convergence can be 
mitigated to some extent by an early exit strategy for the IPM iteration defined
in \eqref{eqn:varuipm}, allowing the iterations to terminate if 
$\norm{\sdiv \uu^{n,l}}/\norm{\sdiv \uu^{n,l-1}} > 1/2$, if this occurs before 
a tolerance or a maximum number of iterations is reached. Of course this only 
improves the slow convergence of IPM with a small inf-sup constant and relaxes the
requirement of choosing an IPM tolerance that is not too small; the entire algorithm
will still ultimately diverge if the IPM does not converge to $\uu^{n,l}$ with a 
sufficiently small divergence. We include a further discussion of the convergence
of IPM for Stokes in subsection \ref{subsec:ipmstokes}.

The problem of non-convex corners causes an issue for solving the transport problem. 
Figure \ref{fig:wspikes} shows the second component of the terminal approximation to 
$\ww$ via \eqref{eqn:varw} on three meshes with comparable mesh sizes. 
The figure on the left is computed on a mesh generated by {\tt mshr} with parameter 
$M=16$. In the center
is the same computation performed on the computational mesh shown in figure 
\ref{fig:cduct-mesh} with two uniform refinements.
On the right of figure \ref{fig:wspikes}
is a variant in which the re-entrant corner is smoothed into two corners with smaller 
internal angles.
The height of the spike is ten times higher on the the mesh generated by
{\tt mshr} than either of the other two. Simulations run on either uniform or boundary
refinements of this mesh realized poor convergence of the IPM iteration, indicating 
degeneration of the inf-sup constant.  
The smoothed domain on the right has multiple spikes which are
moderately shorter than on the mesh shown in \ref{fig:cduct-mesh}.
The singular behavior of $\ww$ at nonconvex corners is a feature of the problem
formulation given by \eqref{eqn:diffrfsgeetoo}. 
For this problem, a good computational mesh should allow for local refinement and
resolution of the singularity without causing the algorithm to diverge.

\begin{figure}
\includegraphics[trim = 100pt 0pt 100pt 0pt,clip = true,width=0.30\textwidth]
{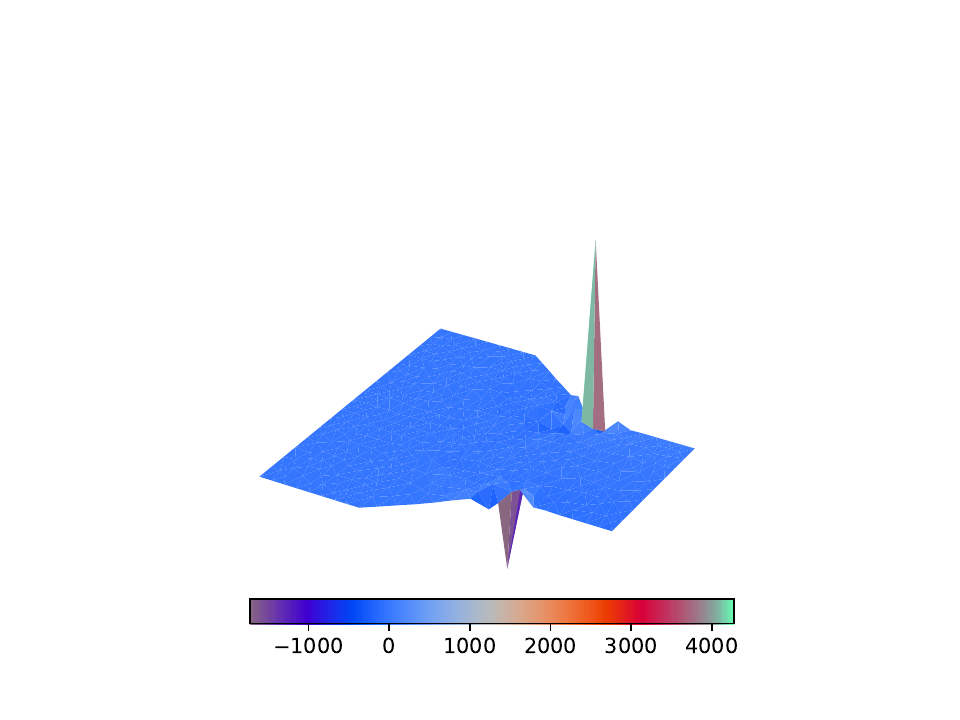}
\includegraphics[trim = 100pt 0pt 100pt 0pt,clip = true,width=0.30\textwidth]
{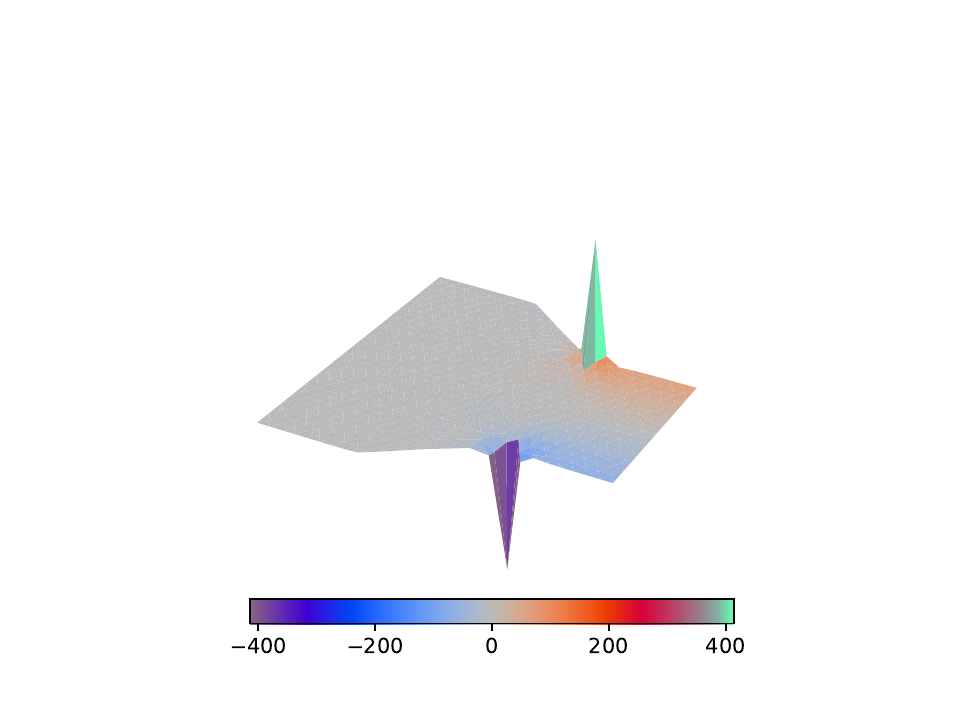}
\includegraphics[trim = 100pt 0pt 100pt 0pt,clip = true,width=0.30\textwidth]
{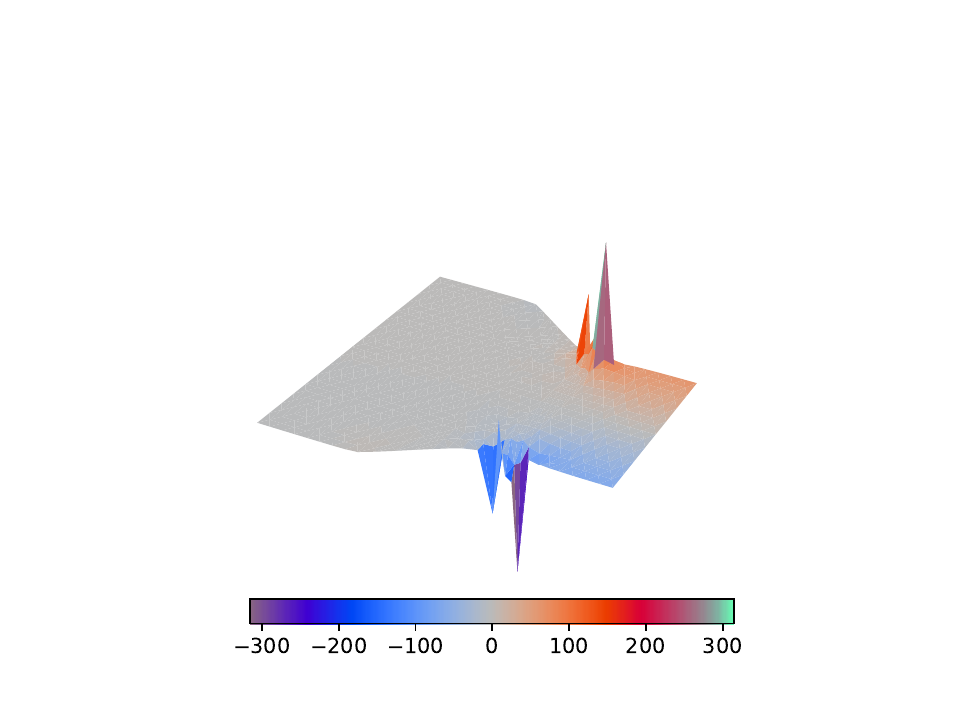}
\caption{Spikes in $w$ at the re-entrant corners of the domain on three different
meshes with $U = 1$, $\nu = 1$, and $\alpha_1 = \alpha_2 = 0.1$.}
\label{fig:wspikes}
\end{figure}

\subsection{IPM and Stokes}\label{subsec:ipmstokes}
We noted in the previous section that the grade-two simulations for Poiseuille flow
may degrade slightly as the mesh is refined and as the iterated-penalty method (IPM)
tolerance is decreased. This is unfortunately a feature of the IPM which may be 
mitigated by defining a mesh on which the inf-sup constant is robust.
\begin{table}
\begin{center} 
\begin{tabular}{|c|c|c|c|c|c|}
\hline
$M$&$\norm{\uu-\uu_h}_{H^1}$ &$\norm{\sdiv \uu_h}_{L^2}$ &IPM iters &$\rho$ & split type \\
\hline
 32 &1.67e-08 &1.02e-11  &   3   &1.00e+04 & crossed-triangle \\
 32 &2.34e-09 &4.90e-11  &   4   &1.00e+03 & crossed-triangle \\
 32 &3.68e-10 &3.91e-11  &   7   &1.00e+02 & crossed-triangle \\
\hline
 32 &4.27e-09 &1.21e-11  &   3   &1.00e+04 & right-triangle \\
 32 &5.76e-10 &5.52e-11  &   4   &1.00e+03 & right-triangle \\
 32 &3.22e-10 &5.63e-11  &   9   &1.00e+02 & right-triangle \\
\hline
\end{tabular}
\end{center}
\vspace{-5mm}
\caption{Stokes errors for Poiseuille flow in the domain \eqref{eqn:chanlosgeetoo} with $L=1$
for $\nu=1$, using the iterated penalty method \cite{lrsBIBih} 
(IPM) with quartics on an $M\times M$ array of squares split in two ways.
The Malkus crossed-triangle split consists of
squares divided into four triangles by the bisectors.
The right-triangle split consists of
squares divided into two right triangles.
}
\label{tabl:justokerho}
\end{table}
Table \ref{tabl:justokerho} indicates the effect of the penalty parameter $\rho$ on
the convergence of IPM.
Larger $\rho$ gives faster convergence for the divergence, but at the
expense of less accuracy for the velocity.
Table \ref{tabl:justokerho} shows the errors after the optimal number of IPM iterations
in terms of minimizing the velocity error, $\uu-\uu_h$.
We see that a smaller error can be achieved at the cost of doing more iterations
with a smaller $\rho$.
As $\rho$ is decreased further, the number of IPM iterations becomes
prohibitively large.

\section{Conclusion}\label{sec:conc} 

We have demonstrated that it is possible, with suitable care, to simulate the
grade-two model in a geometry related to a contraction rheometer.
We have indicated how the results can be used to 
determine the viscosity $\nu$ from experimental data.
We have also explored issues related to identifying the grade-two parameters
$\alpha_1$ and $\alpha_2$ with certain caveats.
In particular, the force data appears to be the same for distinct values of $\balpha$,
however, we identify a regime in which $\balpha$ may be identifiable with a contraction
rheometer.


\section{Acknowlengments}

We thank Mats Stading for valuable discussions.
SP was supported in part by the National Science Foundation NSF DMS-2011519.

\bibliographystyle{plain}
\bibliography{numax}

\appendix

\section{Spaces}
\label{sec:spaces}

Here we collect the notation used for various Sobolev spaces and norms.
We denote by $L^p(\Omega)$ the Lebesgue spaces \cite{lrsBIBgd} of $p$-th
power integrable functions, with norm
$$
\norm{f}_{L^p(\Omega)}=\bigg(\intox{|f(\xx)|^p}\bigg)^{1/p}.
$$
Note that we can easily apply the same notation to vector or tensor valued $f$.
We think of tensors of any arity as vectors of the appropriate length, and we
think of $|f(\xx)|$ as the Euclidean length of this vector.
For tensors of arity 2 (i.e., matrices) this is the same as the Frobenius norm.
We will write the spaces for such tensor-valued functions as $L^p(\Omega)^m$
for the appropriate $m$ (e.g., $m=d^2$ for arity 2).
Similarly, we denote by $L^\infty(\Omega)$ the Lebesgue space of essentially
bounded functions, with
$$
\norm{f}_{L^\infty(\Omega)}=\sup\set{|f(\xx)|}{\hbox{a.e.}\;\xx\in\Omega}.
$$
Correspondingly, we define Sobolev spaces and norms of order $m$ by
$$
\norm{f}_{W^m_p(\Omega)}=\bigg(\sum_{|\alpha|\leq m}\norm{D^\alpha f}_{L^p(\Omega)}^p\bigg)^{1/p},
$$
where $D^\alpha$ is the weak derivative $\partial^\alpha/\partial\xx^{|\alpha|}$
\cite{lrsBIBgd}.
More precisely, the spaces $W^m_p(\Omega)$ are defined as the subspaces of $L^p(\Omega)$
for which the corresponding norm is finite.
The case $p=2$ is denoted by $H$:
$$
H^m(\Omega)=W^m_2(\Omega).
$$

\section{Determining inflow boundary conditions}
\label{sec:detinflo}

The proposed method \eqref{eqn:algfrfsgeetoo} requires
specification of boundary conditions for $\ww=-\Delta\uu+\nabla\pi$.
Here we compute $\ww$ for typical flow geometries.
It corresponds to the divergence of the stress.

\subsection{Grade-two channel flow}
\label{sec:gtwochanflo}

In \cite{lrsBIBjd}, simple two-dimensional grade-two flows (Couette and Poiseuille)
are presented for the domain $\Omega$ defined by
\begin{equation}\label{eqn:chanlosgeetoo}
\Omega=\set{\xx\in\Rtwo}{0<x_1<L,\; 0<x_2<1}.
\end{equation}
Suppose that $u_2\equiv 0$ and $u_1$ depends only on $x_2$.
This is true for shear flow (Couette flow) and pressure-driven flow (Poiseuille flow).
For the remainder of this subsection, we refer to $u_1$ as just $u$
to simplify notation.
For shear (Couette) flow, 
$\ww=\bfz$.
For Poiseuille flow, in the channel \eqref{eqn:chanlosgeetoo}, 
$$
\gbc=\uu= U\begin{pmatrix}x_2(L-x_2)\\ 0 \end{pmatrix} ,\qquad 
\ww =-\frac{2U^2}{\nu}(L-2x_2)\begin{pmatrix} 0 \\ 2\alpha_2+3\alpha_1\end{pmatrix}.
$$
Furthermore,
\begin{equation}\label{eqn:definepeeg}
p(\xx)=-2U\nu  x_1 + (2\alpha_1+ \alpha_2) U^2(L-2x_2)^2 +c_p.
\end{equation}

\subsection{Grade-two pipe flow}

Consider Poiseuille flow in a circular pipe.
To be specific, we define the domain $\Omega$ to be
\begin{equation}\label{eqn:pipelosgeetoo}
\Omega=\set{\xx\in\Rtre}{x_1^2+x_2^2<1,\; 0<x_3<L}.
\end{equation}
Suppose that $u_1=u_2\equiv 0$ and 
\begin{equation}\label{eqn:pipeprofiltoo}
u_3(\xx)=U \big(1- \big(x_1^2+x_2^2\big)\big).
\end{equation}
This is true for pressure-driven flow (Poiseuille flow).
For such flows, $\uu\cdot\nabla\uu=\bfz$, and
the strain rate $\nabla\uu$ is given by
$$
\nabla\uu=-2U\begin{pmatrix}0& 0 & 0\\ 0 & 0 & 0 \\ x_1 & x_2 & 0 \end{pmatrix},\qquad
\nabla\uu^t=-2U\begin{pmatrix}0& 0 & x_1\\ 0 & 0 & x_2 \\ 0 & 0 & 0 \end{pmatrix}.
$$
Thus
$$
\A =-2U\begin{pmatrix}0& 0 & x_1\\ 0 & 0 & x_2 \\  x_1 & x_2 & 0 \end{pmatrix},\qquad
\uu\cdot\nabla\A=\bfz, \qquad
\A\tc\A=
4U^2\begin{pmatrix}x_1^2 & x_1 x_2 &0 \\ x_1 x_2 & x_2^2 &0\\0&0& x_1^2+x_2^2\end{pmatrix},
$$
$$
\A\tc(\nabla\uu)=
4U^2\begin{pmatrix}x_1^2 & x_1 x_2 &0 \\ x_1 x_2 & x_2^2 &0\\0&0& 0\end{pmatrix},\qquad
\A\tc(\nabla\uu^t)=
4U^2\begin{pmatrix}0 & 0  &0 \\ 0 & 0 &0\\0&0& x_1^2+x_2^2\end{pmatrix},
$$
$$
(\nabla\uu)^t\tc\A=\big(\A^t \tc(\nabla\uu)\big)^t=\big(\A\tc(\nabla\uu)\big)^t
=\A\tc(\nabla\uu)=
4U^2\begin{pmatrix}x_1^2 & x_1 x_2 &0 \\ x_1 x_2 & x_2^2 &0\\0&0& 0\end{pmatrix},
$$
$$
\A\tc(\nabla\uu)+(\nabla\uu)^t\tc\A=
8U^2\begin{pmatrix}x_1^2 & x_1 x_2 &0 \\ x_1 x_2 & x_2^2 &0\\0&0& 0\end{pmatrix},
$$
$$
(\nabla\uu)\tc\A=\big(\A^t \tc(\nabla\uu)^t\big)^t=\big(\A\tc(\nabla\uu)^t\big)^t
=\A\tc(\nabla\uu)^t=
4U^2\begin{pmatrix}0 & 0  &0 \\ 0 & 0 &0\\0&0& x_1^2+x_2^2\end{pmatrix}.
$$
We can simplify this by introducing two matrices
$$
\J=
4U^2\begin{pmatrix}x_1^2 & x_1 x_2 &0 \\ x_1 x_2 & x_2^2 &0\\0&0& 0\end{pmatrix},\qquad
\K=
4U^2\begin{pmatrix}0 & 0  &0 \\ 0 & 0 &0\\0&0& x_1^2+x_2^2\end{pmatrix}.
$$
For example,
$$
\A\tc\A=\J+\K,\qquad
\A\tc(\nabla\uu)+(\nabla\uu)^t\tc\A=2\J,\qquad
(\nabla\uu)^t\tc\A=\J.
$$
For the steady-state, grade-two fluid model, the stress tensor simplifies \cite{lrsBIBjd} to
\begin{equation} \label{eqn:tredeeteelook}
\begin{split}
\T_\gs &= \nu \A + \alpha_1\big(\uu\cdot\nabla\A+\A\tc(\nabla\uu)+(\nabla\uu)^t\tc\A\big) 
        + \alpha_2 \A\tc\A\\
&=  \T_\ns
+2\alpha_1 \J +\alpha_2(\J+\K).\\
\end{split}
\end{equation}
The tensor $\btau$ is given  by \eqref{eqn:sigmadefoo}:
\begin{equation}\label{eqn:bisigmaeval}
\begin{split}
\btau &= \alpha_1 (\nabla\uu)^t\tc\A +(\alpha_1+ \alpha_2) \A\tc\A   -\uu\otimes\uu \\
 &= \alpha_1 \J +(\alpha_1+ \alpha_2) (\J + \K)  -\uu\otimes\uu .\\
\end{split}
\end{equation}
We can compute $\sdiv\btau$ as follows.
By definition,
$$
(\sdiv\J)_i=\sum_j J_{ij,j}=J_{i1,1}+J_{i2,2},\qquad
\sdiv\J=\begin{pmatrix} J_{11,1}+J_{12,2}\\J_{21,1}+J_{22,2}\\0\end{pmatrix}
=4U^2\begin{pmatrix} 3 x_1  \\ 3 x_2 \\0\end{pmatrix}
$$
since $\J$ is constant in $x_3$ and thus $\J_{i3,3}=0$.
Similarly, $\sdiv\K=\bfz$.
%
Therefore
$$ 
\sdiv\btau= (2\alpha_1+ \alpha_2) \sdiv \J 
=12U^2 (2\alpha_1+ \alpha_2) \begin{pmatrix} x_1 & x_2 & 0\end{pmatrix}^t.
$$
From \eqref{eqn:divennetoo} we find
\begin{equation}\label{eqn:divennigoobis}
\sdiv N(\uu,\pi)= -\alpha_1 \nabla\uu^t\nabla\pi +\sdiv\btau=
\big(2\alpha_1 U\pi_{x_3} +12U^2 (2\alpha_1+ \alpha_2) \big)
\begin{pmatrix} x_1 & x_2 & 0\end{pmatrix}^t.
\end{equation}
Also
$
-\nu\Delta\uu=4\nu U \begin{pmatrix} 0 & 0 & 1\end{pmatrix}^t.
$
Thus
$$
4\nu U \begin{pmatrix} 0 \\ 0 \\ 1\end{pmatrix}
+\nabla p = \sdiv\big(\T_\gs-T_\ns\big)=\sdiv\big((2\alpha_1+\alpha_2)\J)\big)
=12U^2 (2\alpha_1+ \alpha_2) \begin{pmatrix} x_1 \\ x_2 \\ 0\end{pmatrix}.
$$
Therefore
$$
p_{x_3}=-4\nu U,\qquad \begin{pmatrix} p_{x_1} \\ p_{x_2}\end{pmatrix}=
12U^2 (2\alpha_1+ \alpha_2) \begin{pmatrix} x_1 \\ x_2 \end{pmatrix}.
$$
These equations are solved by
$$
p(\xx)= -4\nu U x_3 + 6 U^2 (2\alpha_1+ \alpha_2) \big(x_1^2+x_2^2\big).
$$
Let us make the ansatz that $\pi(\xx)=-4 U x_3 + f(x_1,x_2)$.
We have
\begin{equation}\label{eqn:pwhatwrite}
p=\nu\pi+\alpha_1 u \pi_{x_3}=p- 6 U^2 (2\alpha_1+ \alpha_2) \big(x_1^2+x_2^2\big)
+\nu f - 4U \alpha_1 u \, ,
\end{equation}
which suggests that
\begin{equation}\label{eqn:duweneedis}
\begin{split}
\nu f &=6 U^2 (2\alpha_1+ \alpha_2) \big(x_1^2+x_2^2\big)
+ 4U^2 \alpha_1 \big(1-\big(x_1^2+x_2^2\big) \big) \\
 &= U^2 (8\alpha_1+ 6\alpha_2) \big(x_1^2+x_2^2\big) + 4U^2 \alpha_1 .
\end{split}
\end{equation}
Therefore
$$
\pi(\xx)=-4 U x_3 + \nu^{-1} U^2\big((8\alpha_1+ 6\alpha_2) \big(x_1^2+x_2^2\big)
+ 4 \alpha_1 \big).
$$
In particular, \eqnn{eqn:divennigoobis} implies
\begin{equation}\label{eqn:divennagin}
\sdiv N(\uu,\pi)= U^2 \big(16\alpha_1+12 \alpha_2 \big)
\begin{pmatrix} x_1 & x_2 & 0\end{pmatrix}^t.
\end{equation}
Similarly,
$$
\ww=-\Delta\uu+\nabla\pi=\nu^{-1} U^2(16\alpha_1+ 12\alpha_2)
 \begin{pmatrix} x_1 \\ x_2 \\ 0\end{pmatrix}.
$$
This provides an appropriate initial condition for solving the third 
(transport) equation in \eqref{eqn:algfrfsgeetoo} for pipe flow
with boundary conditions given by \eqref{eqn:pipeprofiltoo}.

Note that $\ww_{,x}\equiv 0$.
Thus \eqref{eqn:diffrfsgeetoo} implies that
$$
\sdiv N = \nu \ww,
$$
which is consistent with \eqref{eqn:divennagin}.


\end{document}